\documentclass[nalpha-refs]{oup-contemporary}

\usepackage{graphicx}
\usepackage{siunitx}

\usepackage{subcaption}
\usepackage{amsmath,amssymb}
\usepackage{nomencl}
\usepackage{longtable,tabularx,float}
\usepackage{hyperref,array,makecell,multirow,booktabs}
\renewcommand{\vec}[1]{\boldsymbol{#1}}
\DeclareMathOperator*{\argmax}{argmax} 

\makenomenclature
\newlength{\nomgroupstartsep}
\usepackage{etoolbox}
\renewcommand\nomgroup[1]{%
  \itemsep\nomgroupstartsep%
  \item[\bfseries
  \ifstrequal{#1}{A}{Alphabets}{%
  \ifstrequal{#1}{B}{Greek Letters}{%
  \ifstrequal{#1}{C}{Dimensionless Groups}{}}}%
  ]
  \itemsep\nomitemsep
}

\journal{jcde}
\title{Airfoil Optimization using Design-by-Morphing}
\doi{10.1093/jcde/xxxxxxx}
\jyear{202X}
\jvolume{X}
\jnumber{X}
\jstartpage{X}
\jendpage{X}
\jpubdate{Day Month Year}

\author[1,{\authfn{1},{\authfn{2}}}]{Haris M. Sheikh}
\author[1,{\authfn{2}}]{Sangjoon Lee}
\author[1]{Jinge Wang}
\author[1]{Philip S. Marcus}

\affil[1]{Department of Mechanical Engineering, University of California, Berkeley, CA 94720, USA}

\authnote{\authfn{1}Corresponding author. E-mail: \href{mailto: harissheikh@berkeley.edu}{harissheikh@berkeley.edu}}
\authnote{\authfn{2}Equal contribution.}

\papercat{Research article}

\runningtitle{Airfoil Optimization using Design-by-Morphing}

\begin{document}
\begin{frontmatter}
\maketitle
\begin{abstract}
Design-by-Morphing (DbM) is a novel design methodology that creates a search space for topology optimization. Traditional design techniques often impose geometric constraints and, sometimes, the designer's biases on the design space, which restricts the novelty of the designs and allows for only small local changes. On contrary, we show in this paper that DbM does not impose such restrictions on the design space, thus allowing for a radical and expansive search space with only a few design parameters. We compare DbM with other methods in the case of design space generation for 2D airfoils and found that DbM can reconstruct the entire UIUC database to >99.5\% accuracy. Furthermore, using a bi-objective genetic algorithm, we optimize the airfoil designs created by DbM to maximize both the lift-over-drag ratio, $CLD_{max}$, and stall angle tolerance, $\Delta \alpha$, which results in a Pareto-front of innovative airfoils that exhibit substantial improvements in both objectives.


\end{abstract}

\end{frontmatter}


\setcounter{section}{0}
\section{Introduction}
Airfoil shape optimization is a critical stage in the design of aerodynamic components, such as aircraft wings \parencite{wing_1, wing_2, wing_3, wing_4} and wind-turbine blades \parencite{wind_turbine_1,wind_turbine_2,wind_turbine_3,wind_turbine_4,wind_turbine_5,wind_turbine_6}. The airfoil optimization process typically involves three main components: shape parameterization, airfoil evaluation, and optimization. Among these, the parameterization method defines both the design space and the complexity of the optimization problem. To ensure effectiveness, a desirable parameterization technique must be able to encompass a wide design space using a modest amount of design parameters \parencite{sobester_barrett_2008,sripawadkul_comparison_2010,masters_review_2015,10.1115/1.4036134}. This is particularly important during the initial design phase, where minimum geometric constraints are imposed, and the flexibility to make significant changes during optimization is beneficial.

Shape parameterization methods differ in their fidelity and control ranges \parencite{masters_review_2015,sobester_barrett_2008}, and can be placed on a \textit{virtual} spectrum according to the geometric scope of each design parameter. At one end of the spectrum, adjusting a single parameter alters a \textit{local} section of the airfoil, which offers precise shape control but modifies the shape slowly. At the opposite end, each design parameter affects the \textit{global} contour of the airfoil \parencite{sobester_barrett_2008}.

At the \textit{local} end of the spectrum is the \text{discrete method} \parencite{discrete_method}, where the design parameters are exactly the discrete points that define the airfoil surface. Since the position of each point can be adjusted, the design space is potentially limitless \parencite{survey_parameterization}, and precise local control with high fidelity can be achieved. 
However, a substantial number of surface points are needed to accurately describe an airfoil shape, which complicates the optimization problem. Gradient-based optimizers are frequently employed to mitigate the increased complexity, but they are likely to get stuck at a sub-optimal solution during the optimization.

As the geometric scope of each parameter is expanded, there emerge the classical approaches that are based on the curve-fittings of regional features or control points.
For example, the popular \text{parametric section (PARSEC) method} \parencite{PARSEC} uses eleven or twelve parameters to represent major sectional features of an airfoil, including leading edge radii and upper and lower crest locations, and constructs the airfoil surface using a $6^{\mathrm{th}}$ order polynomial. Another popular method is the \text{B{\'e}zier parameterization} \parencite{bezier}, which constructs the upper and the lower surfaces of the airfoil through the B{\'e}zier curves defined by pre-chosen control points. Additionally, a hybrid of the two techniques, \text{B{\'e}zier-PARSEC parameterization}, was introduced by \textcite{bezier-parametric}, which uses the parameters of the PARSEC method to define the B{\'e}zier curves that form the shape contours. One main issue with the above methods is their inability or inefficiency to include high-fidelity features; the PARSEC and the B{\'e}zier-PARSEC methods both have a fixed number of parameters and limited range of fidelity, while the B{\'e}zier parameterization requires higher-degree B{\'e}zier curves to describe complex shapes which are inefficient to calculate \parencite{survey_parameterization}. 

To consider finer details of airfoils or, equivalently, to represent more complex curves, either \text{B-splines} \parencite{sanaye_multi-objective_2014, han_adaptive_2014} or \text{nonuniform rational B-spline (NURBS)} \parencite{nurbs} can be used, which creates curves by connecting low-order B{\'e}zier segments defined by control points.
As the number of control points increases, these methods move to the local end of the spectrum and become capable of representing high-fidelity features, but the computing complexity also increases. One way to reduce the number of the design parameters is to group the control points together so that the global transformations such as twisting and thickening can be used as the parameters. This is known as the \text{free-form deformation (FFD) method} \parencite{FFD1,FFD2} and is closer to the \textit{global} end of the spectrum. A similar method, called \text{radial basis function domain element (RBF) approach} \parencite{RBF1,RBF2,10.1115/1.4046650}, also exists and makes use of radial basis function to exert deformation on the airfoil.

Near the \textit{global} end of the spectrum, we see methods using spectral construction of basis functions to form or deform airfoil shapes. One popular choice of the basis functions is the dominant modes from \text{singular value decomposition (SVD)} of an airfoil dataset \parencite{POD1,POD2,PCA, doi:10.2514/6.2020-2707, doi:10.2514/6.2019-1701}. Other choices include sinusoidal functions of the Hicks-Henne approach \parencite{hicks_wing_1978}, which create `bumps' on a reference airfoil surface, and surface functions of the class/shape function transformation (CST) method \parencite{kulfan_fundamental_2006, akram_cfd_2021}, which are in the form of the product of a class function and a shape function generated by a linear combination of Bernstein polynomials. Nonetheless, like many other methods on the spectrum, these methods also suffer from the so-called {\it curse of dimensionality} that more basis functions or modes are always required to resemble high-fidelity features.

Efforts have been made to overcome the {\it curse of dimensionality} 
\parencite{viswanath_forrester_keane_2011, viswanath_forrester_keane_2014, cinquegrana_iuliano_2018}. A recent work by \textcite{chen_chiu_fuge_2020} applied a generative adversarial network (GAN) to learn the major shape variations of an airfoil database and use those to parameterize the shapes while also preserving the high-fidelity features via an additional noise space. However, like many other dimension reduction methods, this study assumes that the optimum design is not far from the database, which is not always true. To address this limitation, they proposed another GAN-based method that encourages diversity during sample generation \parencite{10.1115/1.4048626}, but a large dataset is still required to initialize the training. In contrast, our paper is motivated by the optimization problem during the early design stage when few initial designs are available. Therefore, we are interested in a parameterization method that is capable of representing high-fidelity features even when the design parameters and initial airfoil designs are limited.

In this paper, we apply the Design-by-Morphing (DbM) parameterization technique to the airfoil optimization problem. DbM is a novel and universal design strategy that was first introduced by \textcite{2018CompM6223O} and has been used in recent years for geometry optimization of different problems \parencite{2018CompM6223O,2019APS..DFDQ14007S,sheikh_2022,2021APS..DFDA15004S}. As a \textit{global} method, it `morphs' homeomorphic baseline shapes together to create new shapes and is able to interpolate and extrapolate the design space, allowing for both high-fidelity representation of shapes without the curse of dimensionality and radical modifications to the shapes without any implicit geometric constraints \parencite{2018CompM6223O,sheikh_2022}. This strategy is applicable to a variety of 2D and 3D design problems and we aim to conduct a special case study of DbM for the 2D airfoil shape optimization here. Throughout this paper, we aim to make the following scientific contributions:
\begin{itemize}[noitemsep,topsep=0pt]
    \item Application of DbM to 2D airfoil shape optimization, showing its accurate reconstruction of the existing airfoil database and radical changes in airfoil shapes while being free from geometric constraints and designers' biases by extrapolation of the design space by applying negative weights.
    \item Evaluation of airfoil design capacity of the DbM strategy and comparison with other typical 2D airfoil design strategies.
    \item Sensitivity analysis for the number of baseline shapes, convergence analysis compared to conventional airfoil design strategies is shown and the significance of extrapolation for DbM is also shown.
    \item Optimization within the 2D airfoil search space generated by DbM using a genetic algorithm and investigation of the optimum Pareto-front.
\end{itemize}

\section{Design-by-Morphing}\label{sec:DbM}
Design-by-Morphing (DbM) works by morphing homeomorphic, i.e. topologically equivalent, shapes to create a continuous and constraint-free design search space. 
It comes with several advantages. To begin with, DbM is valid for shapes of any dimensions, and capable of creating exotic shapes because radically different baseline shapes can be morphed together. Furthermore, DbM does not impose any geometric constraints on the design parameters. And the only implicit constraints are the selections of the `baseline shapes' themselves, which are necessary to prescribe the problem to be solved. Lastly, it is able to create an extensive design search space, even when the number of pre-existing designs is small, e.g. \parencite{sheikh_2022}, by both the means of `extrapolation', that is to assign negative weights during morphing, and the inclusion of irregular or uncommon shapes. The details of DbM for airfoil optimization are presented in the subsequent subsections.

\subsection{Baseline Shapes and Morphing}
The DbM technique requires two or more homeomorphic `baseline shapes', mostly chosen from pre-existing designs in the literature, to create the design space. A one-to-one correspondence between the baseline shapes must first be established through some systematic shape collocation methods in either the functional \parencite{2018CompM6223O} or the geometric space \parencite{2019APS..DFDQ14007S,sheikh_2022}. Then the new shapes can be generated by applying weights to the collocation vectors of the baseline shapes and summing them together in a linear manner. 

\begin{figure}[bt!]
\begin{subfigure}[bt!]{\linewidth} 
\vbox{
\centering{
  \includegraphics[width=.9\linewidth]{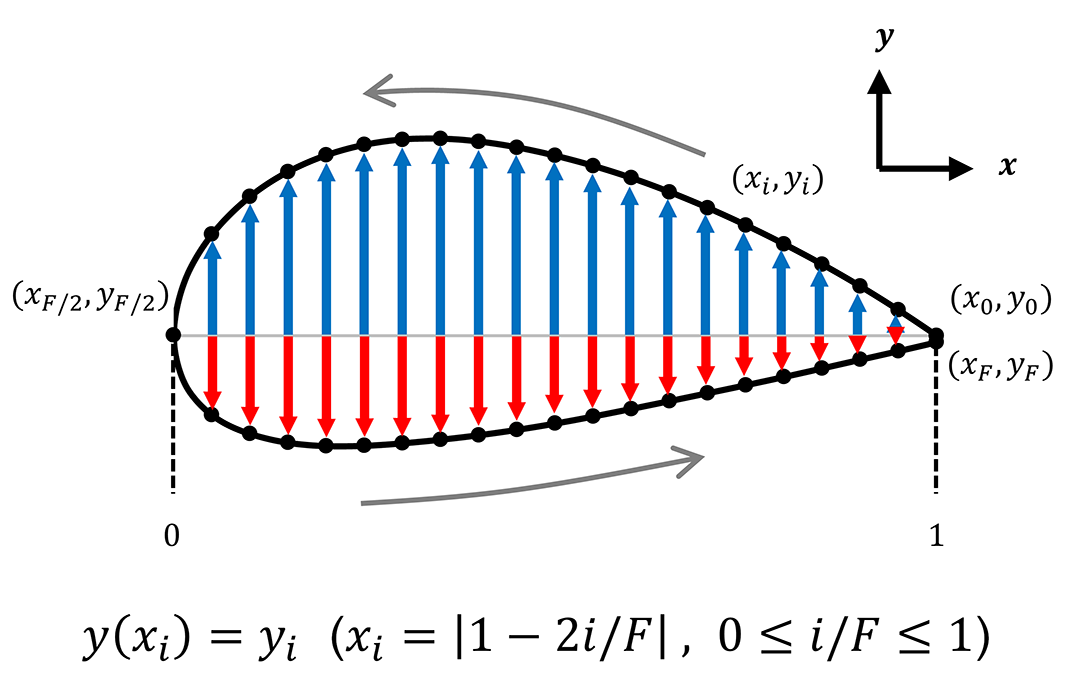}
}%
}%
\subcaption{\centering\label{fig:DbM-coll_a} $n^{\text{th}}$ airfoil shape}
\end{subfigure}
\\[1.7em]
\begin{subfigure}[bt!]{\linewidth} 
\centering{
\includegraphics[width=.9\linewidth]{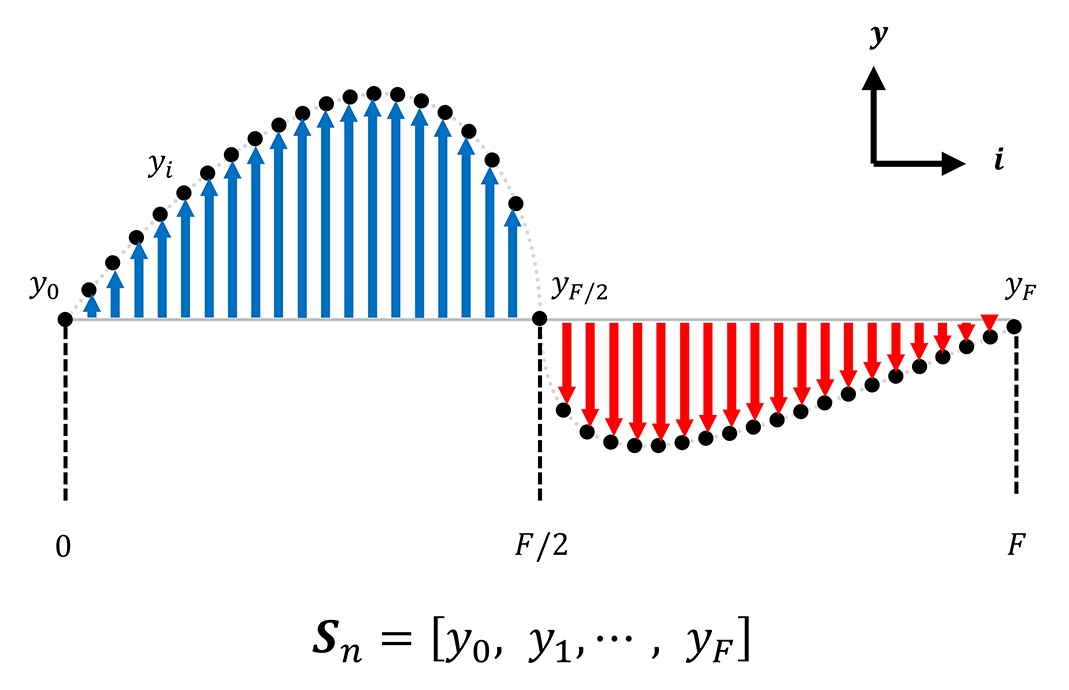}
\subcaption{\centering\label{fig:DbM-coll_b} $y$-coordinate collocation vector}
}\end{subfigure}%
\caption{An example of DbM. The coordinates of the baseline shapes are weighted, summed, and normalized to form the coordinates of a morphed shape.}\label{fig:DbM-coll}
\end{figure}

For 2D airfoils, the closed shapes can be collocated in the Euclidean coordinate system. It is noted here that all 2D shapes bounded by a single surface are homeomorphic to one another. Using the leading edge of each airfoil as origin, each shape can be collocated by taking fixed and uniformly spaced points along the $x$-axis
, creating a one-to-one correspondence between the shapes. This collocation strategy is demonstrated in Figure~\ref{fig:DbM-coll}, and the baseline shapes used in this paper are chosen from various airfoils in the literature, which are detailed later. Morphing is performed by multiplying a specific airfoil shape with a scalar weight, summing the weighted vectors, and then normalizing them. For a collection of $N$ baseline shapes, morphing is given by
\begin{equation}
\vec{P}(\vec{x}) = \dfrac{1}{\sum
_{m=1}^{N}w_m}\sum\limits_{n=1}^{N} w_n \vec{S}_n(\vec{x}) \, .
\label{eq:1}
\end{equation}

\begin{figure*}
    \centering
    \includegraphics[width=.9\linewidth]{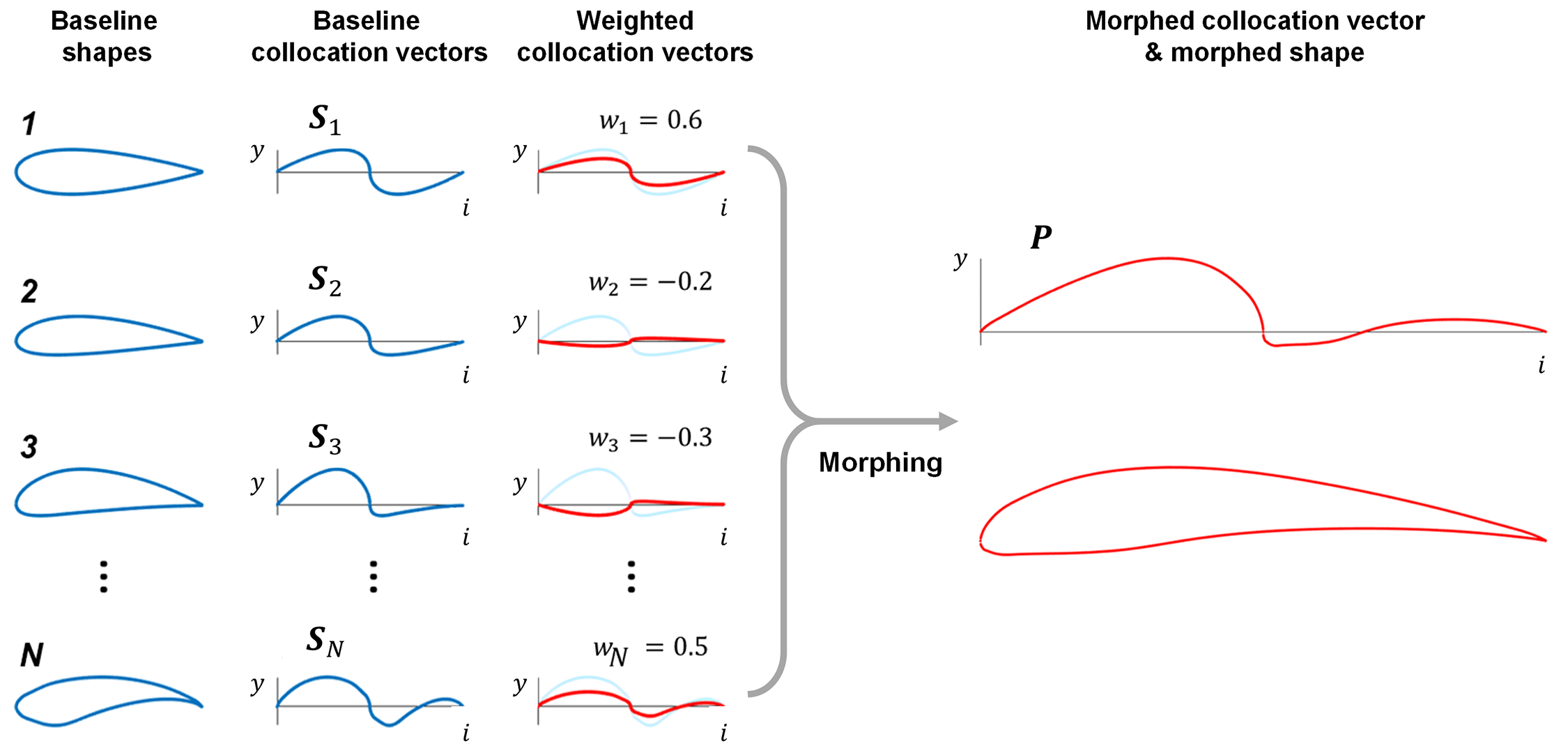}
    \caption{Application of DbM to 2D airfoils. Column 1 shows the baseline shapes. Column 2 depicts the elements of the collocation vectors of the baseline shapes plotted as a function of the index $i$ of the collocation vector. Column 3 shows the weighted elements of the collocation vector plotted as a function of the index $i$ of the collocation vector. Column 4 shows the resultant collocation vector of the morphed shape and the morphed shape itself.}
    \label{fig:DbM-strat}
\end{figure*}

\noindent Here $\vec{S}_n(\vec{x})$ is the $y$-coordinate collocation vector of the $n^{\text{th}}$ baseline shape, collocated at $\vec{x} = \left[ x_0, \cdots , x_F \right]$ where the $i^\text{th}$ $x$-coordinate $x_i = \left| 1-2i/F \right|$ and $F$ is the number of collocation points. The first half of the elements of $\vec{S}_n$ represents the top surface of the airfoil, and the second half of the elements of $\vec{S}_n$ represents the bottom. $w_n \in [-1,1]$ is the morphing weight applied to the $y$-coordinate vector of the $n^{\text{th}}$ baseline shape, and negative $w_n$ values imply extrapolation. A visual demonstration of the strategy is presented in Figure~\ref{fig:DbM-strat}. 

\subsection{Intersection Control}
\begin{figure}[tb!]\centering{
\begin{subfigure}[t]{0.495\linewidth}
\vbox{
\centering{
  \includegraphics[width=1.015\linewidth]{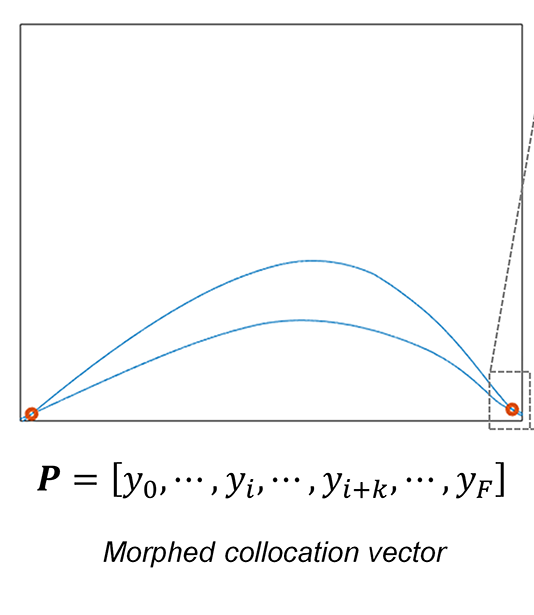}
}%
}%
\subcaption{\centering\label{fig:DbM-intersect_a}}
\end{subfigure}
\begin{subfigure}[t]{0.495\linewidth}
\centering{
\includegraphics[width=1.015\linewidth]{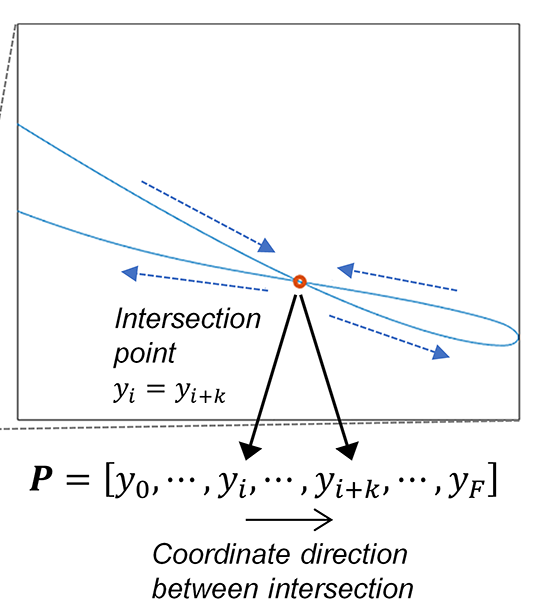}
\subcaption{\centering\label{fig:DbM-intersect_b}}
}\end{subfigure}\\[0.5em]%
\begin{subfigure}[t]{0.495\linewidth}
\vbox{
\vspace*{1em}%
\centering{
  \includegraphics[width=1.015\linewidth]{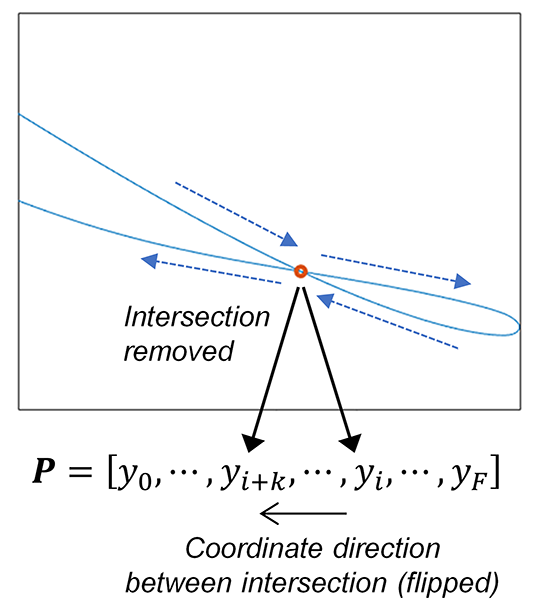}
}%
}%
\subcaption{\centering\label{fig:DbM-intersect_c}}
\end{subfigure}
\begin{subfigure}[t]{0.495\linewidth}
\centering{
\includegraphics[width=1.015\linewidth]{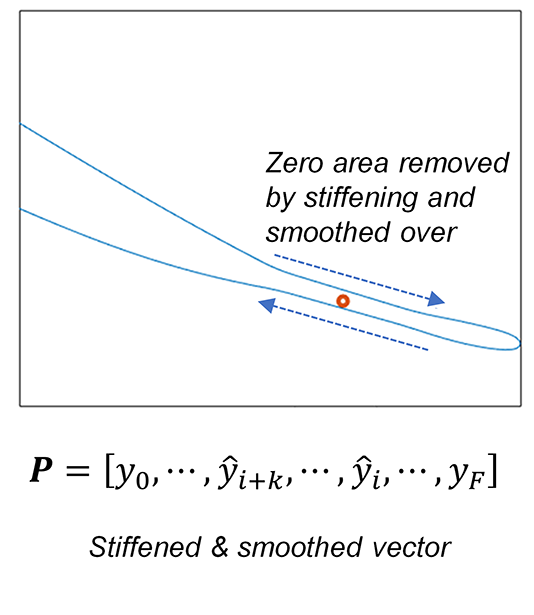}
\subcaption{\centering\label{fig:DbM-intersect_d}}
}\end{subfigure}%
\caption{Conditioning for intersection removal. (\subref{fig:DbM-intersect_a}) Intersections are detected; (\subref{fig:DbM-intersect_b}) A blown-up image of one intersection, with the shape coordinates direction depicted by arrows; (\subref{fig:DbM-intersect_c}) Intersection removed by flipping the vector between the intersection; (\subref{fig:DbM-intersect_d}) The `zero-area' removed by linear interpolation and then smoothed over, as shown by hatted y-coordinates.}
\label{fig:DbM-intersect}}
\end{figure}

For smooth baseline shapes, applying positive weights, i.e. interpolation, will always create smooth shapes. However, applying negative weights, i.e. extrapolation, may produce non-physical geometries, such as self-intersections, which have `zero-area' regions, as shown in Figure~\ref{fig:DbM-intersect_a}. One may discard the morphed airfoil shapes with self-intersections during the optimization, but that diminishes the size of our design space. Instead, we recover new shapes by removing the intersections. 

Intersection removal is accomplished by first locating the intersection within the morphed coordinate vector and restructuring the vector by `flipping' it between the intersection points, as shown in Figure~\ref{fig:DbM-intersect_c}. The vector is then `stiffened' to remove the `zero-area' between the intersections by removing the points in their neighborhoods and linearly interpolating between the broken coordinate vectors. As seen in Figure~\ref{fig:DbM-intersect_d}, this removes the `zero-area' space and adds some physical area to the shape at the original intersection point. This process is repeated until all intersections are removed, e.g. both intersections in Figure~\ref{fig:DbM-intersect} are successfully removed, and finally, a moving-average smoothing filter is applied to smooth out any sharp edges.

\subsection{Selection of Baseline Shapes}
The selection of baseline shapes is a crucial component of the DbM strategy and ultimately determines the size and novelty of our search space. Metaphorically, the selection of the baseline airfoil shapes serves as the gene pool for the morphed airfoils and its diversity is important for creating a large design space. 

One way of selecting the baseline shapes is by performing singular value decomposition (SVD) or principle component analysis (PCA) on a set of shapes and then using the dominant modes as the baseline shapes. Methods such as parametric model embedding \cite{SERANI2023115776} can help reduce the dimensionality of the problem as well. Although these methods would help in quantitatively choosing baselines, these methods, however, require an existing dataset which might not be available in many shape optimization problems \parencite{2018CompM6223O, doi:10.1177/0309524X17709732, sheikh_2022}. Therefore, while techniques like SVD and PCA can be easily applied to airfoil shape optimization problems, and provide arguably better baselines, we choose the baseline shapes qualitatively instead to demonstrate the universality of DbM even for engineering problems with few existing designs. In other words, for research purposes we assume that the airfoil database is not \textit{a priori} knowledge at the selection stage, except for those chosen as baselines.

An additional benefit of directly morphing existing designs is that, from a human designer's perspective, it can be more intuitively informative than handling PCA modes. For example, vertical-axis wind turbines are broadly categorized into drag, lift and hybrid categories, so the weights associated with each type are more informative to a human designer than the weights of the dominant modes. On the other hand, choosing actual shapes as the baseline shapes has the advantage that non-conforming designs can be easily added, as is the case for baseline \#19 (\textit{mirrored} Selig airfoil). Conventional techniques may have much more difficulty adding radical features into the design space, and the significance of such radical baseline shapes is demonstrated in the Results section. 

\begin{figure}[tb!]
    \centering
    \includegraphics[width=\linewidth]{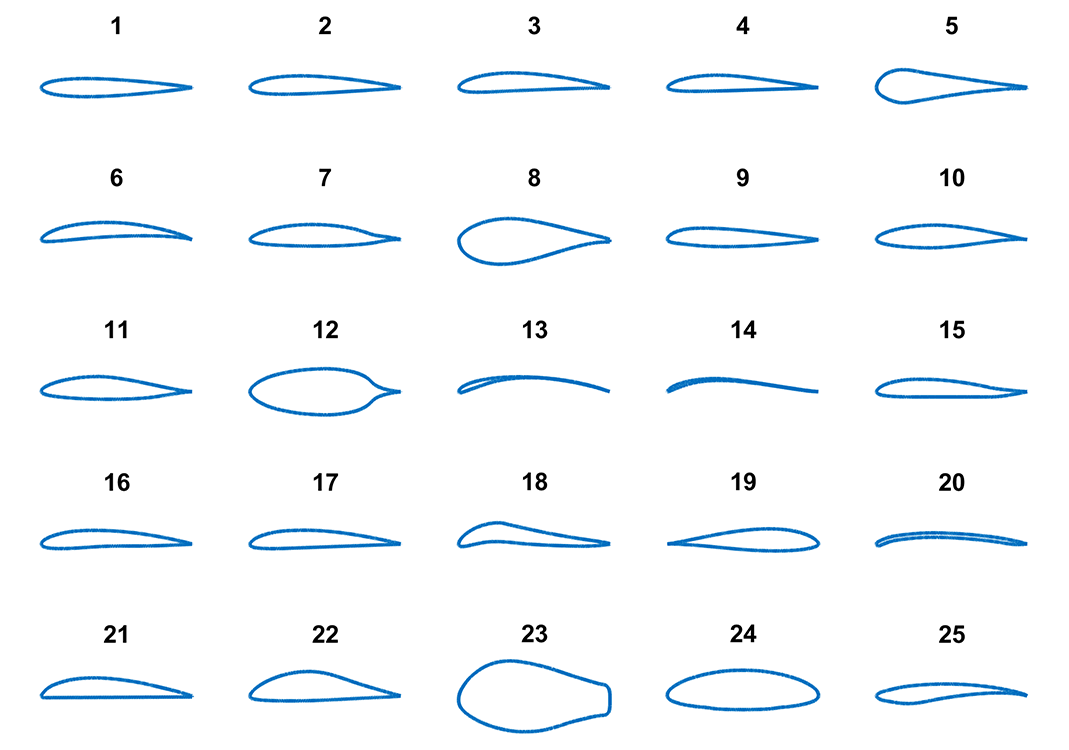}
    \captionof{figure}{Twenty-five baseline airfoil shapes chosen from the UIUC database \parencite{uiuc_database}. See Appendix \ref{app:baseline} for further details.}
    \label{fig:baseline}
\end{figure}

We selected 25 shapes (see Figure~\ref{fig:baseline}) from the UIUC airfoil coordinates database \parencite{uiuc_database} as our baseline shapes. They were picked to ensure diversity and to introduce radical features into the design space. Our selection of baseline shape included airfoils that are either known for high lift-to-drag ratio or good stall performance, which are commonly used in the literature and the industry. We also included airfoils with poor aerodynamic performances, as well as airfoils with irregular shapes, to provide novelty to the design space. It is worth noting that, unlike in the conventional airfoil optimization processes \parencite{koroglu_ozkol_2019}, we deliberately included the bad performers so that our optimization could suppress these features by assigning them negative weights, which will be demonstrated in greater detail in our later results. The model names and characteristics of the baseline shapes can be found in Appendix \ref{app:baseline}. To express shapes as collocation vectors, each airfoil shape is represented by 4,001 coordinates that span counterclockwise from the upper surface trailing edge around the leading edge to the lower surface trailing edge, with equally distributed $x$-coordinates parallel to the airfoil chord line of unit length (i.e. $F=4,000$ in Figure~\ref{fig:DbM-coll}).


\subsection{Airfoil Design Capacity Test}
As a benchmark, we reconstruct the entire UIUC airfoil database \parencite{uiuc_database} using DbM to test the robustness of our method and the representation capacity of the generated design space. Having noted that one of the key features of DbM is to permit shape extrapolation, we compare our reconstruction results against the results of an interpolation-only DbM (DbM-I) where all DbM weights are non-negative. In addition, we performed the same test on 3 conventional 2D airfoil shape parameterization methods: PARSEC \parencite{PARSEC}, NURBS \parencite{nurbs} and the Hicks-Henne approach \parencite{hicks_wing_1978}. These tests were meant to answer the following questions:
\begin{itemize}[noitemsep,topsep=0pt]
    \item How much does the extrapolation expand the design space?
\end{itemize}
Extrapolation is undoubtedly better at creating a wider design search space. However, we shall focus on how \textit{quantitatively} the search space is broadened by extrapolation, so as to confirm whether this feature genuinely distinguishes DbM from other generic approaches.
\begin{itemize}[noitemsep,topsep=0pt]
    \item Is DbM comparable to conventional airfoil shape parameterization methods in terms of shape reconstruction?
\end{itemize}
It should be noted that we selected the baseline airfoil shapes for DbM solely based on the qualitative principle of ensuring diversity and intentionally avoided the use of a selection method that requires a known, rich design database in advance. We shall demonstrate that the answer to the above question is still positive, even though DbM is not a design method specifically for airfoils, unlike the methods that the DbM is compared with.

\begin{figure}[tb!]
    \centering
    \includegraphics[width=\linewidth]{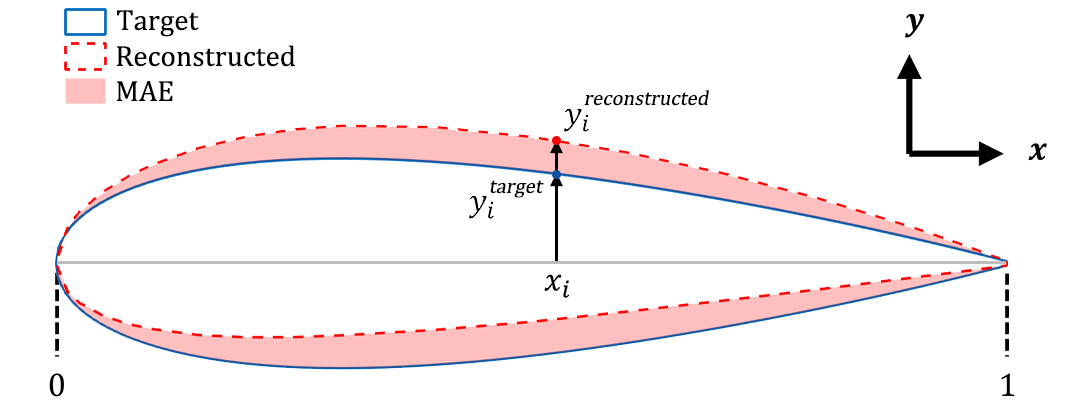}
    \captionof{figure}{Geometric demonstration of mean absolute error (MAE) between target and reconstructed airfoil surfaces.}
    \label{fig:mae_example}
\end{figure}

For all of the 1,620 airfoils in the UIUC database \parencite{uiuc_database}, we obtained the closest representation of each target shape by running a global optimization of the input design parameters that minimizes mean absolute error (MAE) between the target and the reconstructed airfoil surfaces.
A geometric demonstration of MAE is provided in Figure~\ref{fig:mae_example}. Using the functional expression of an airfoil $y(x_i)= y_i$ ($0\le i \le F$), as introduced in Figure ~\ref{fig:DbM-coll_a},
\begin{equation}
\text{MAE}(target,~reconstructed) = \frac{2}{F} \int_0^F{y_i^{{error}} di} ,
\label{eq:xor}
\end{equation}
where $y^{{error}}_i \equiv | y^{{reconstructed}}_i - y^{{target}}_i|$. 
When we express the error in percentage terms, i.e. (MAE $\times$ 100) \%, we emphasize that the error is described as a proportion of the area difference to the square of the chord length, as all airfoil shapes are normalized to maintain a unit chord length. The factor of $2$ in Eq. \eqref{eq:xor} is present for this reason.

\begin{table*}[]
\caption{Airfoil shape parameterization methods for comparison}\label{table:comparison_methods}
\centering{%
  \begin{tabularx}{\linewidth}{c C C c}
    \toprule 
   \makecell[c]{Method} & \makecell[c]{Design Variables (DVs)} & \makecell[c]{\# of DVs} &  \makecell[c]{Remark} \\\midrule
    PARSEC & $\begin{aligned}  r_{le}^{up/lo}:~& \text{Leading edge radii} \\[-0.5ex] x^{up/lo},~ y^{up/lo}:~& \text{Crest coordinates}\\[-0.5ex] y_{xx}^{up/lo}:~& \text{Crest curvatures}\\[-0.5ex] y_{te},~ t_{te}:~& \text{Trailing edge mid-position and thickness}\\[-0.5ex] \alpha_{te},~ \beta_{te}:~& \text{Trailing edge direction and wedge}  \end{aligned}$ & 12 & \makecell[c]{Fixed \# of parameters} \\ \\
    
    NURBS  & $\begin{aligned} x_{ctrl,i}^{up/lo}, y_{ctrl,i}^{up/lo}:~& \text{Control point coordinates}~(i=1,\cdots,~4) \\[-0.5ex] w_{i}^{up/lo}:~& \text{Curve weights}~(i=1,\cdots,~4) \\[-0.5ex] y_{te}^{up/lo}:~& \text{Trailing edge positions} \end{aligned}$  & 26 & \makecell[c]{Third-order B-spline \\ Evenly distributed knots} \\ \\
    
    Hicks-Henne & $\begin{aligned} w_{i}^{up/lo}:~& \text{Bump widths}~(i=1,\cdots,~6) \\[-0.5ex] m_{i}^{up/lo}:~& \text{Bump magnitudes}~(i=1,\cdots,~6) \end{aligned}$ & 24 & \makecell[c]{Base profile: NACA 0012 \\ Cosine-distributed bump points} \\ \\
    
    DbM (Present) & $\begin{aligned} w_{i}:~& \text{Morphing weights}~(i=1,\cdots,~25) \end{aligned}$ & 25 & \makecell[c]{See Figure~\ref{fig:baseline} for the baselines} \\[-1.5ex] \\\bottomrule 
  \end{tabularx}
  }%
  \end{table*}

To obtain the closest representation, we utilized a MATLAB-based single-objective genetic algorithm (GA): \texttt{ga}. The population size is set to 100, and the maximum number of generations is set to 500. The lower bound for MAE was set to $1.44 \times 10^{-3}$ (or 0.144 \%) from  Eq. \eqref{eq:xor} for a chord length of 1, in accordance with the lower limit of Kulfan's typical wind-tunnel tolerance \parencite{kulfan_fundamental_2006,Masters2017}. To ensure a fair comparison, all the airfoil parameterization methods tested underwent the same optimization scheme for shape reconstruction, with similar numbers of design parameters (e.g. 25 design parameters for the DbM) except for PARSEC, which has fixed design parameters. In general, the fidelity of these design methods improves as the number of design variables increases \parencite{Masters2017}.

It is important to note that the objective of our reconstruction tests is to examine both the efficiency and the accuracy of a given parameterization method in the context of shape generation during the design process, which must be distinguished from the accuracy in surface fitting. NURBS, for instance, can achieve arbitrary accuracy for shape fitting if a good initial guess of the parameters is provided, but it may not be ideal for shape generation as the shape it constructs varies slowly during the optimization process. Accordingly, all the reconstructions tests are initialized in a consistent manner to provide a meaningful comparison. In particular, the initial population is set to contain a single parameter set that represents the profile of NACA 0012 with the remaining sets randomly distributed. The results of the reconstruction tests can then be understood as the ability of a method to create various shapes, including the common designs that have been collected in the UIUC database, precisely within a certain number of optimization generations. And the progressive improvement in the design space reconstruction can be observed as a function of number of GA generations.


\begin{figure}[tb!]
    \centering{%
    \includegraphics[width=\linewidth]{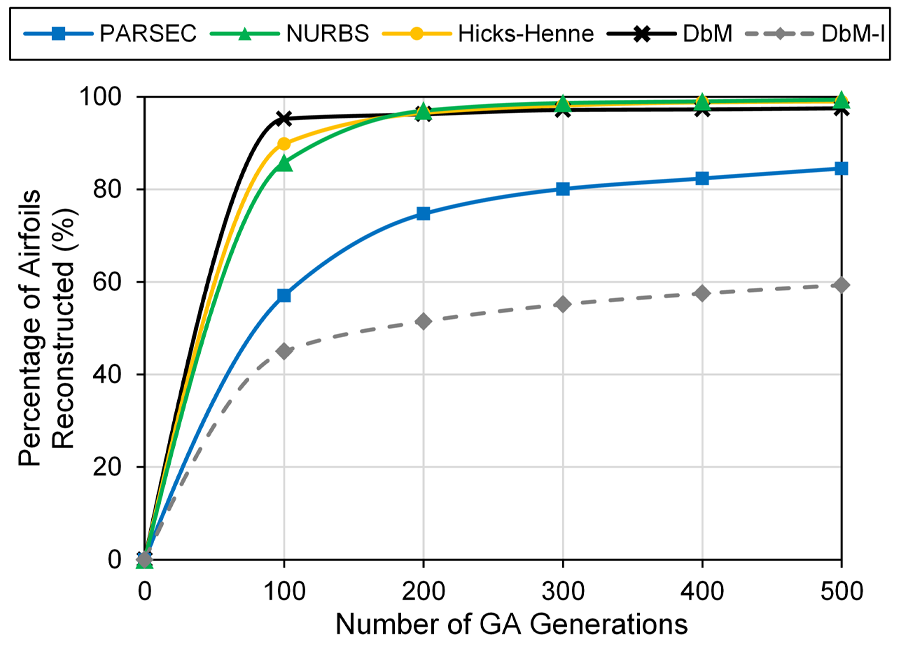}
    \caption{Percentage of airfoils reconstructed within the mean absolute error tolerance of 0.5 \% using DbM, compared to those from DbM only with interpolation (DbM-I) and three airfoil shape parameterization methods (see Table~\ref{table:comparison_methods}).}\label{fig:comparison} 
}%
\end{figure}

Applying DbM for the reconstruction of the UIUC database, we found that 1,618 of the 1,620 airfoils of the entire UIUC database, were reconstructed with a MAE error <~1\%. Even for the two airfoils with the highest error, the DbM reconstruction still resulted in a MAE error of less than 1.5 \%. Figure~\ref{fig:comparison} displays the percentage of airfoils that were reconstructed within the tolerance of 0.5 \% MAE error, with respect to the total number of GA generations. A comparison between DbM and DbM-I reveals that the extrapolation feature of DbM significantly contributes to the improved performance of the method, suggesting that the extrapolation feature is indispensable for DbM. On the other hand, at the maximum GA generation, the total percentage of reconstructed airfoils with a MAE error <~0.5\% increases from 60 \% (DbM-I) to 98 \% (DbM). As a result, DbM is converges faster than any other conventional approaches tested here.

\begin{figure*}\centering
\begin{subfigure}[t]{0.195\linewidth}
\vbox{
\centering{
  \includegraphics[width=\linewidth]{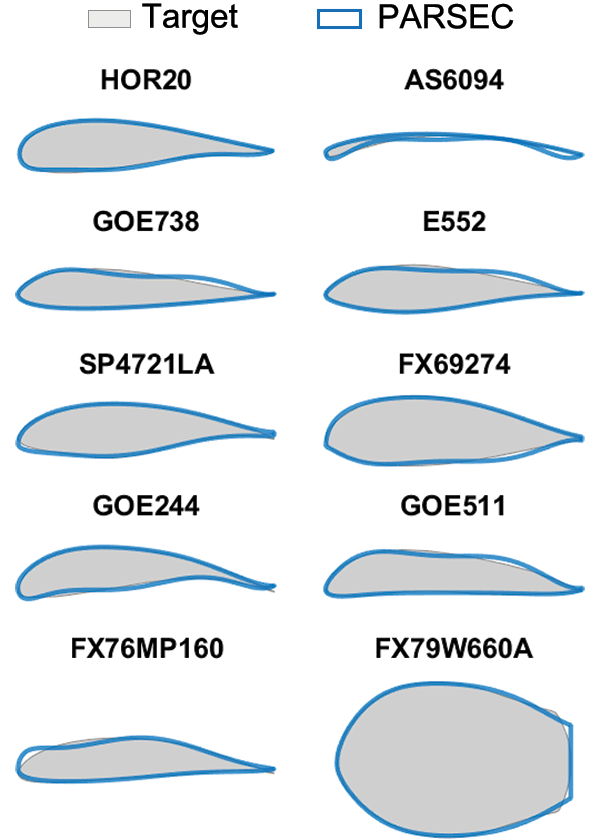}
}%
}%
\subcaption{\centering\label{fig:DbM-reconstruct_a} PARSEC}
\end{subfigure}
\begin{subfigure}[t]{0.195\linewidth}
\centering{
\includegraphics[width=\linewidth]{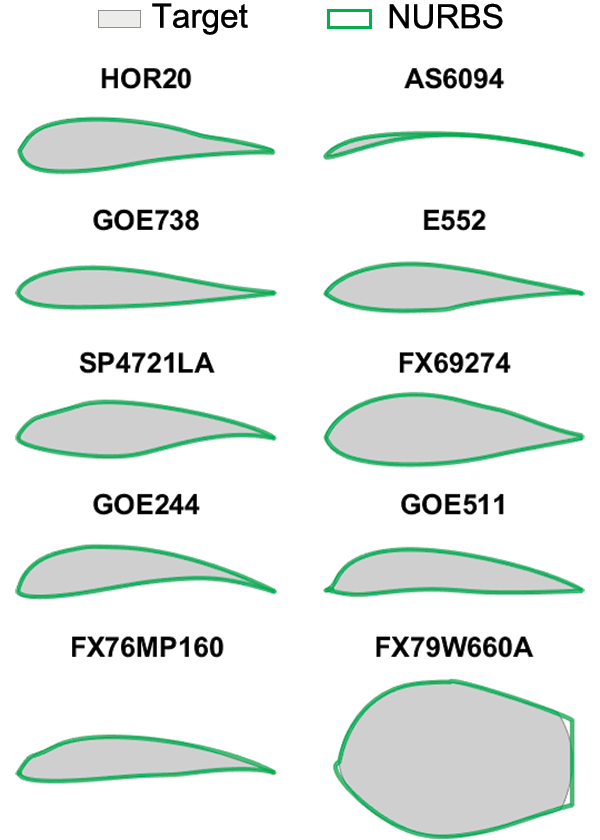}
\subcaption{\centering\label{fig:DbM-reconstruct_b} NURBS}
}\end{subfigure}%
\begin{subfigure}[t]{0.195\linewidth}
\vbox{
\centering{
\includegraphics[width=\linewidth]{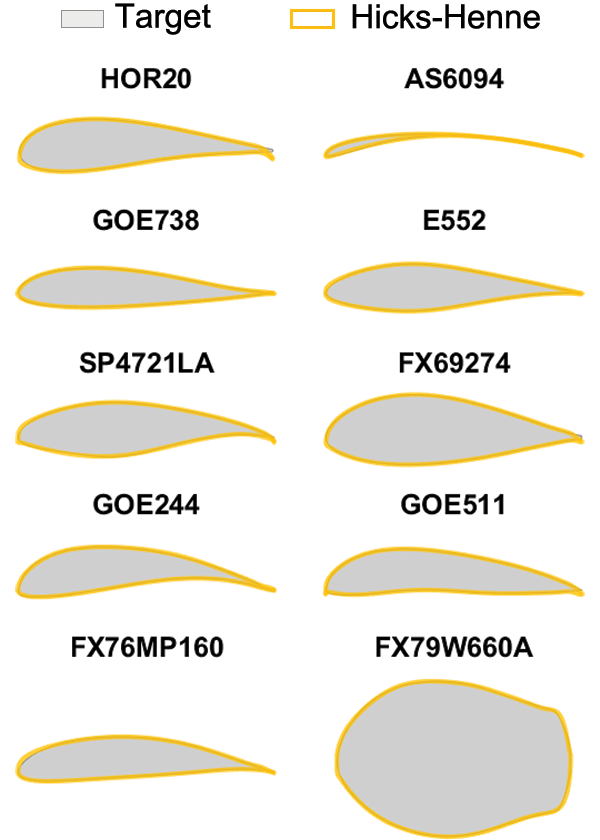}
}%
}
\subcaption{\centering\label{fig:DbM-reconstruct_c} Hicks-Henne}
\end{subfigure}
\begin{subfigure}[t]{0.195\linewidth}
\vbox{
\centering{
\includegraphics[width=\linewidth]{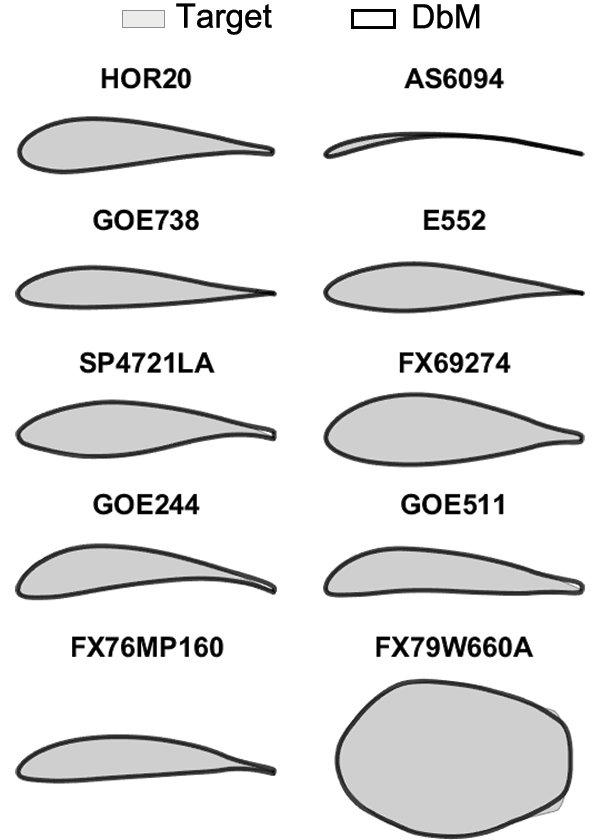}
}%
}
\subcaption{\centering\label{fig:DbM-reconstruct_d} DbM}
\end{subfigure}
\begin{subfigure}[t]{0.195\linewidth}
\vbox{
\centering{
\includegraphics[width=\linewidth]{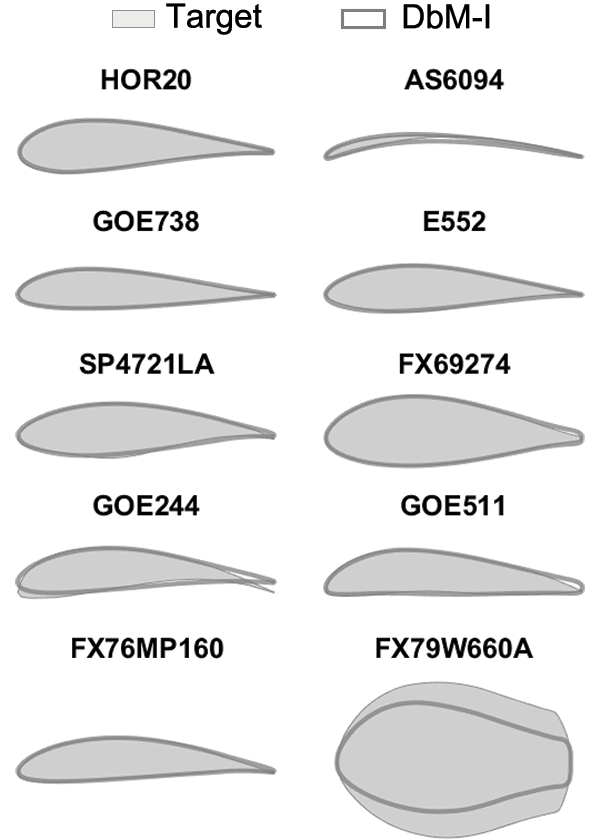}
}%
}
\subcaption{\centering\label{fig:DbM-reconstruct_e} DbM-I}
\end{subfigure}
\caption{GA-based reconstruction of pre-existing airfoil shapes using different design parameterization methods, after 500 GA generations in total.}\label{fig:DbM-reconstruct}
\end{figure*}

To provide better insight, we plot the reconstruction results of 10 airfoil shapes in Figure~\ref{fig:DbM-reconstruct}, which represent the less successful reconstruction attempts. In particular, we ranked the results based on the unweighted average of all MAE errors from the five tested methods for each airfoil case and made the selections at every percentile bin from the worst. These 10 shapes are depicted in row-major order, from one in the $90^{\mathrm{th}}$ percentile (HOR 20, average MAE 0.35 \%) to one in the $99^{\mathrm{th}}$ percentile (FX 79-W-660A, average MAE 1.1 \%). Even these less successful results appear to reasonably reconstruct the target airfoil shapes. It is worth paying attention to the worst case, FX 79-W-660A, which is designed for use on a thick rotor blade of a wind turbine and far from the typical streamlined airfoil shapes. DbM-I encountered a notable failure in this case because none of our chosen baseline airfoil shapes for DbM were as thick as the target shape. As a result, the reconstruction just ended up with the thickest baseline airfoil shape, \#23 (see Figure~\ref{fig:DbM-reconstruct_e}). This specific example underscores the significance of the extrapolation feature of DbM, which provides the opportunity to explore extraordinary designs, such as much thicker airfoils in this case. We also observed that NURBS occasionally produced thorn-like local structures (e.g. at the leading edge of GOE 511 in Figure~\ref{fig:DbM-reconstruct_b}), which resulted from the locally deforming nature of NURBS. These artifacts are normally removed by fitting softwares, such as \texttt{FitCrv} in Rhinoceros 3D, or manual handling of the control points by the designer.



\begin{figure}[tb!]
    \centering
    \includegraphics[width=\linewidth]{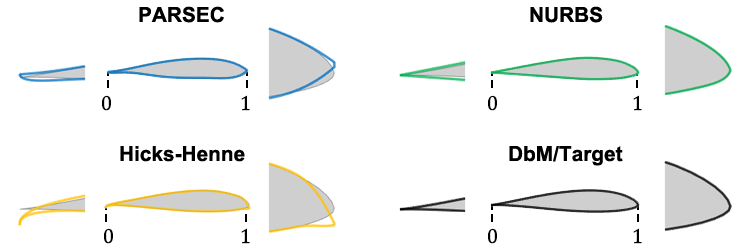}
    \captionof{figure}{Reconstruction of a deliberately ``mirrored" airfoil shape. In contrast to convention, the blunt edge is at $x=1$ and the sharp edge is at $x=0$. PARSEC and the Hicks-Henne method, which implicitly define edge geometries, face challenges in reconstructing the mirrored edges. NURBS looks to perform well due to its better flexibility in adjusting curvatures through weights. DbM has no problem because it can take such non-conforming shapes as baselines as needed.}
    \label{fig:comparison_mirror_case}
\end{figure}

Another way DbM explicitly introduces novelty is by using novel shapes directly as baselines. Generally speaking, novel designs that contain unconventional features can be challenging to construct. For example, the `mirrored' airfoil in our DbM baseline shapes (\#19) is considered off-design by conventional airfoil parameterization methods that prescribe fixed edge geometries for the airfoil, such as relatively `blunt' and `sharp' edges at $x=0$ and $x=1$, respectively. Figure~\ref{fig:comparison_mirror_case} displays the results of GA-based reconstruction of the mirrored airfoil using the tested methods. PARSEC and the Hicks-Henne approach, which implicitly define edge geometries for airfoils, clearly struggle in reconstructing the mirrored edges. At $x=1$, these methods still exhibit hints of the sharp edge in their reconstructed shapes. On the other hand, NURBS performs well as it is more flexible in handling curvatures through weight parameters. However, although not considered in this study, if one wishes to introduce a tentative higher-order feature, such as a stepped wing \parencite{lumsdaine_investigation_1974}, NURBS may require a larger number of design parameters (i.e. more control points and weights). On the contrary, DbM would only require one additional design parameter (i.e. by adding it as a new baseline shape) to introduce novelty regardless of the complexity of the new design.

\begin{figure}[tb!]
    \centering
    \includegraphics[width=\linewidth]{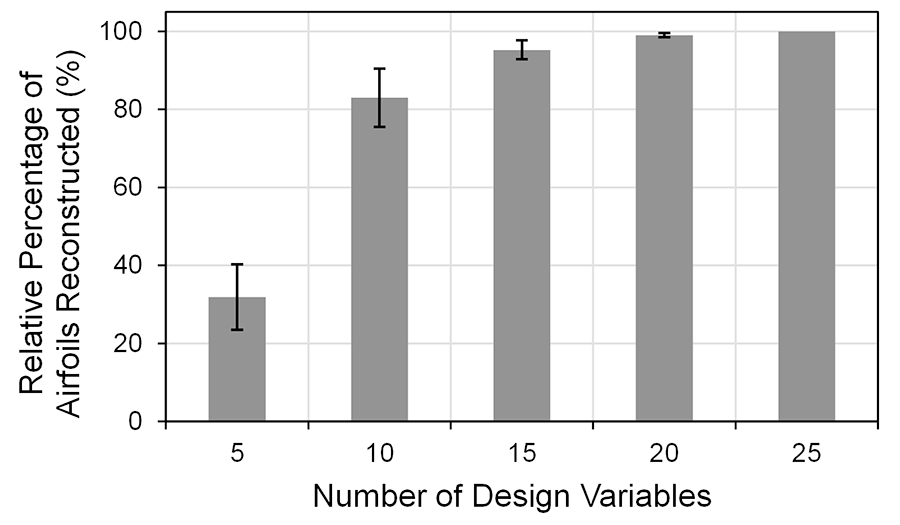}
    \captionof{figure}{Relative percentage of airfoils reconstructed within the mean absolute error tolerance of 0.5 \% by DbM using a subset of the chosen baseline shapes (see Figure~\ref{fig:baseline}) to those using all twenty-five baseline shapes at the maximum GA generation. All tests were done five times with five random subsets for each case. Error bars indicate the standard deviation of the five test results.}
    \label{fig:dimensionality}
\end{figure}

For the current study, we note here that we used 25 baselines based on the computational budget available and our study shows that the number of baseline shapes was sufficient. However, a smaller number of baseline shapes might have proven to be enough as well. to this end, Figure~\ref{fig:dimensionality} presents a sensitivity study of DbM in relation to the number of design variables (baseline shapes) used, where the convergence trend confirms that 20-25 baseline shapes are sufficient. Note that the current sensitivity study is restricted to the same baseline shape set in Figure~\ref{fig:baseline}. All tests were done five times with five random subsets for each case to consider sensitivity to the choice of baseline shapes within the subset. Future efforts will be directed towards conducting additional sensitivity analyses of DbM by varying the selection of the initial 25 baseline shapes themselves.

Overall, we have shown that DbM is competitive against conventional local parameterization methods despite being on the global end of the spectrum. In addition, DbM's ability to generate extraordinary designs through the extrapolation feature enhances the chances of finding novel solutions that deviate from the inputs (baseline shapes), which is important because the aerodynamic performance of an airfoil can be non-intuitively correlated with geometric features of the airfoil. More importantly, DbM is not just a method for design parameterization, but rather a universal design strategy for broader design search. While we have compared DbM to airfoil shape parameterization methods in the context of airfoil optimization, DbM can be useful for any type of problem that aims to introduce more novelty in design search.

\section{Optimization Methodology}\label{sec:method} 

\begin{figure}[tb!]
    \centering
\includegraphics[width=.8\linewidth]{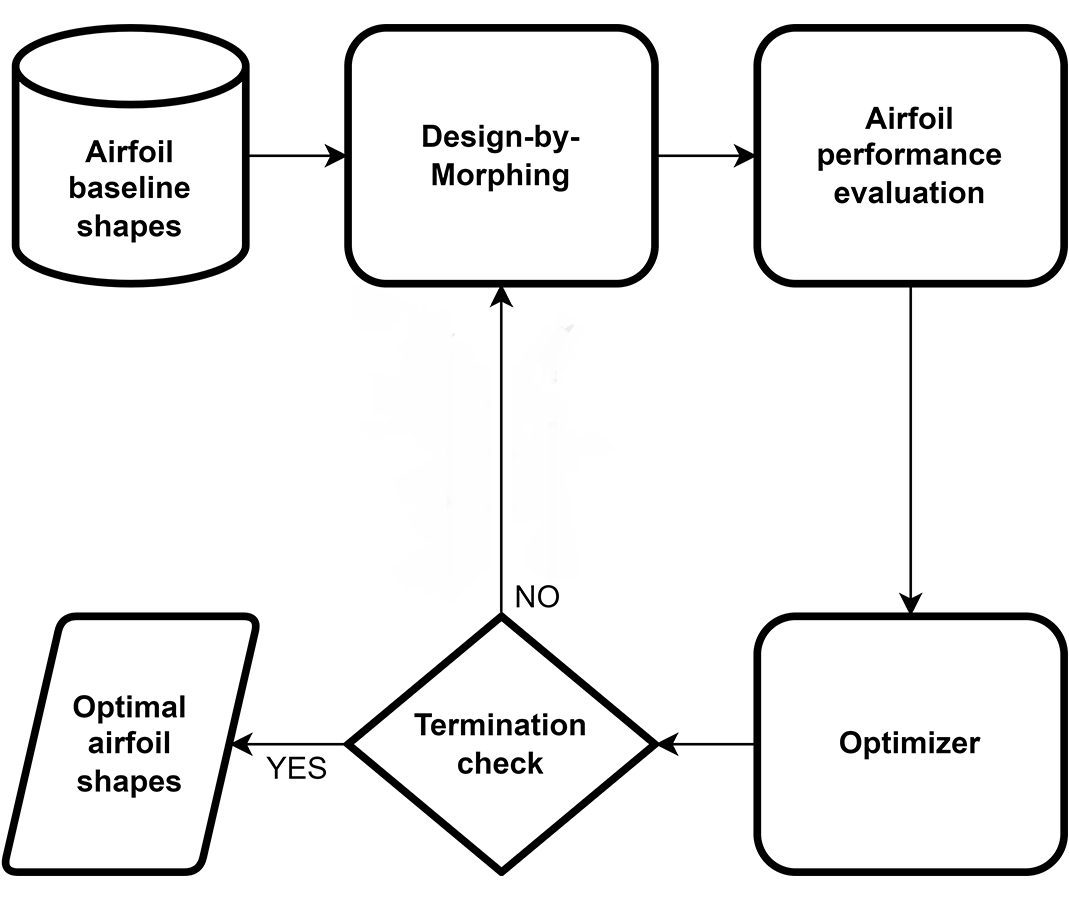}
    \captionof{figure}{General flowchart of airfoil optimization via DbM.}
    \label{fig:flowchart}
\end{figure}

Our airfoil optimization methodology is built around the DbM technique introduced in Sec.~\ref{sec:DbM}. As shown by the flowchart in Figure~\ref{fig:flowchart}, the optimization starts with the selection of the baseline shapes and then evaluates and optimizes the airfoils formed by morphing these baseline shapes using DbM. Our methodology does not rely on a specific airfoil evaluation tool or a specific optimizer, and discussions on their choices are provided in Sec.~\ref{sec:method.evaluation} and Sec.~\ref{sec:method.optimization}, respectively.

\subsection{Airfoil Evaluation} \label{sec:method.evaluation}
Our optimization methodology is not limited to a specific airfoil performance analyzer. Any reliable CFD or experimental methods can be used. For optimizing airfoil shapes using CFD-based solvers, the evaluation of the objective functions (i.e. aerodynamic properties) is typically divided into two categories: the full Reynolds-averaged Navier-Stokes (RANS) based approach and the interacted viscous/inviscid zonal approach. The RANS-based approach is computationally expensive and demands a highly efficient optimizer. To accommodate a large number of design variables, which is common in the aerodynamic designs, a gradient-based optimizer coupled with adjoint methods for computing derivatives is deemed most feasible \parencite{rans_piotrowski_2022,rans_he_2019,rans_kenway_2016}. On the other hand, the viscous/inviscid zonal approach, which combines separated solutions for inviscid external flow and viscous shear layer flow in an iterative manner to form a continuous profile, is faster and less expensive. 

Among a number of inviscid/viscous zonal airfoil analysis codes, XFOIL \parencite{drela_1989} has been the most dominant and widely adopted program \parencite{xfoil_use_ronsten_1992,xfoil_use_selig_1998,xfoil_use_jones_2000,xfoil_use_mueller_2003,xfoil_use_johnson_2005,xfoil_use_batten_2006,xfoil_use_lafountain_2010,xfoil_use_ramanujam_2017,10.1115/1.4048626}. It combines a vorticity panel method for exterior flow with an integral boundary-layer method for viscous boundary layers and uses an $e^9$-type amplification formulation to determine the transition point \parencite{drela_1989}. Its suitability for airfoil designs has been demonstrated in the past literature, where its predictions of aerodynamic properties are found to be in good agreement with both wind-tunnel experiment data \parencite{xfoil_McGhee_1988,xfoil_selig_1995} and the RANS-based simulation results \parencite{xfoil_morgado_2016}. 

While our choice of the evaluation tool is flexible, for this work, we opt for XFOIL due to its acceptable accuracy under our flow condition and its low computation cost. Its widespread usage also allows for a quick reproduction of our optimization results. It is used in a black-box manner so that any other commercial or in-house airfoil analysis tools can be incorporated into our optimization framework if necessary. Our detailed airfoil evaluation setup is given in Appendix \ref{app:baseline}.


\subsection{Optimization} \label{sec:method.optimization}

When a set of solutions is given, the most optimal solution within the set can be determined without difficulty for single-objective optimization problems, which is the case for most of the previous airfoil optimization studies  \parencite{10.1115/1.4046650, 10.1115/1.4054631,10.1115/1.4048626}. However, for multi-objective optimization, multiple and potentially conflicting objectives must be considered simultaneously to determine the optimal answer in the solution set \parencite{Kais99, doi:10.1080/23311916.2018.1502242}. If the designer has a quantitative ranking of the objectives, these objectives can be combined together to formulate a single-objective problem, but when no such ranking exists, constructing a Pareto-front is the most common methodology \parencite{barron, das, mgp}, which has applications in the design of architected materials \parencite{doi:10.1142/S1758825120501161, doi:10.1126/sciadv.abk2218} and turbo-machinery \parencite{inproceedings, 10.3389/fenrg.2021.708230, 2016, RODRIGUES2016587,2011CMAME.200..883W,10.1115/1.2779899}, process engineering \parencite{Nguyen_2021, HuangGAO85, 10.1093/bioinformatics/btz544}, shape design \parencite{LI2020525,app11020729,CIARDIELLO2020115984}, and structural engineering \parencite{10.3389/fmars.2019.00017, article}.

We pose the multi-objective optimization problem as
\begin{equation}
\vec{w}_{opt}= \argmax_{\vec{w} \in \mathcal{W}} (\vec{f}(\vec{w})),
\label{eq:multiobj}
\end{equation}
where $\vec{f}(\vec{w})=[f_1(\vec{w}),f_2(\vec{w}), \cdots ,f_{\text{\textit{K}}}(\vec{w})]$. Here $f_1, \cdots ,f_{\text{\textit{K}}}$ are the $K$ objectives to be maximized, and $\vec{w}$ is the design variable vector. Generally, $\vec{w}$ is a $d$-dimensional vector defined over a bounded set $\mathcal{W}\subset \mathbb{R}^d$ representing $d$ continuous variables. $\{\vec{w}_{opt}\}$ is a set of Pareto-optimal solution vectors, i.e. vectors that are not Pareto-dominated by any other vectors. For the reader's convenience, it is noted that a design variable vector $\vec{\hat{w}}$ is Pareto-dominated by another design variable vector $\vec{\tilde{w}}$ if $f_{k}(\vec{\hat{w}}) \leq f_{k}(\vec{\tilde{w}})$ for all $k \in \left\{ 1, \cdots,K \right\}$. To obtain the Pareto-front, especially when objectives cannot be weighted or when a non-convex black-box function is considered, 
evolutionary or genetic algorithms are a natural choice \parencite{10.1093/bioinformatics/btz544, 860052}. In fact, these algorithms have been commonly implemented in many previous aerodynamic optimization studies due to their gradient-free nature and wide search domain \parencite{akram_cfd_2021, skinner_state---art_2018, rahmad_single-_2020, nsga_zhao_2014}. An alternate choice is Bayesian optimization method, which has been proven to be efficient when the cost functions are expensive to compute (e.g. when using experiments or CFD as an evaluation tool) \parencite{mixmobo}.

Our study considers a bi-objective ($K=2$) two-dimensional airfoil shape optimization. In particular, we optimize the shape of a subsonic airfoil operating in an incompressible flow with $\mathrm{Re} \equiv {U c}/{\nu}$ of $1\times 10^6$, where $U$ and $\nu$ are the free-stream flow speed and fluid kinematic viscosity, respectively, and $c$ is the airfoil chord length. The parameter to be optimized is the morphing weight vector for the DbM technique, defined as:
\begin{equation}
    \vec{w} \equiv (w_1, \cdots , w_{25}) \in \mathcal{D}^{25}, 
\end{equation}
where $\mathcal{D} = [-1,1]\subset \mathbb{R}$ and $w_i$ ($i = 1, 2, \cdots, 25$) is the weight applied to the $i^{\mathrm{th}}$ baseline shape. The design objectives are the maximum lift-drag ratio over all possible angles of attack $\alpha$, i.e. $f_1(\vec{w})=CLD_{max}(\vec{w})$, 
and the difference between the stall angle $\alpha_s$ and the angle where the maximum lift-drag ratio occurs, i.e. $f_2(\vec{w})=\Delta\alpha(\vec{w})$, often called the stall angle tolerance. This particular combination of design objectives has applications in the design of vertical-axis wind turbines \parencite{2019APS..DFDQ14007S}, and the precise definitions of these design objectives are explained in Appendix~\ref{app:prelim}. Both objectives are evaluated using XFOIL, which is efficient enough to be used with the GA.

We use a MATLAB-based variant of the popular NSGA-II algorithm \texttt{gamultiobj} \parencite{996017}, which is a controlled, elitist genetic algorithm. Its practical employment can be found in \textcite{keane_robust_2020} for the purpose of airfoil design optimization, as in the case with ours. Our initial population consists of the single-objective optimums of each design target and random samples in the design space. A population size of 372 is used with a total of 3,000 GA generations. Within each generation, solutions are actively ranked to maintain diversity and to prevent over-crowding in the Pareto-optimal solution set. Our setup was tested using the commonly used set of `$ZDT$' benchmark problems for multi-objective problems, as suggested by \textcite{10.1162/106365600568202}. The test problems and validation results are detailed in Appendix~\ref{app:testfunc}. 

\begin{figure}[tb!]
  \centering
  \includegraphics[width=\linewidth]{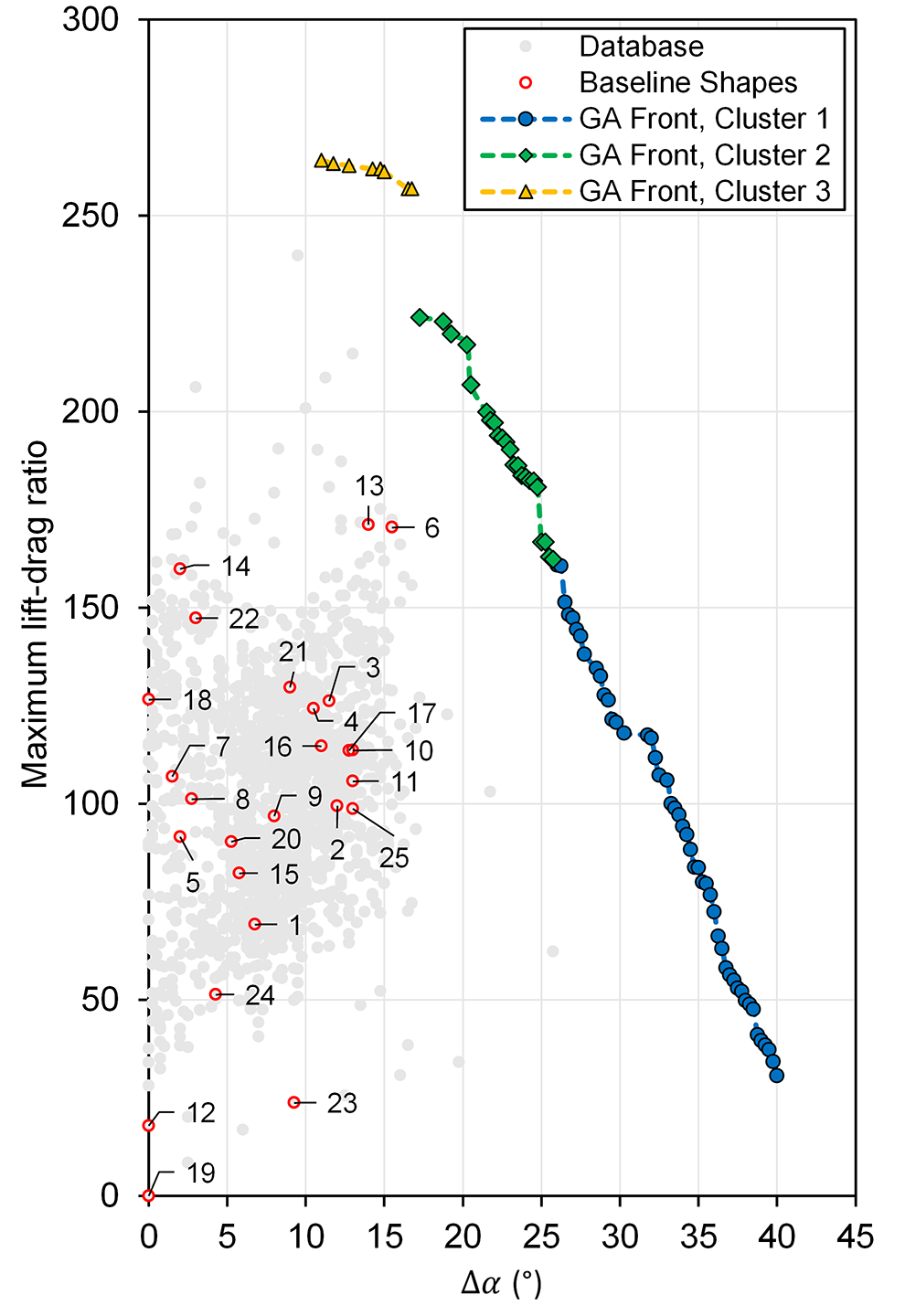}
  \caption{The Pareto-front consisting of the optimal airfoil shapes, resulting from the NSGA-II 3,000 generation runs. The gray points are the whole evaluation outcomes of the UIUC reference database \parencite{uiuc_database}, with the DbM baseline cases in the present study highlighted as red hollow circles with respective indices. See Appendix \ref{app:pca} for the clustering. $\mathrm{Re_{chord}} = 1\times 10^6$}.
  \label{fig:pareto_front}
\end{figure}

\section{Results}
A Pareto-front on the $\Delta \alpha$ - $CLD_{max}$ objective plane resulting from a total of 3,000 generations of the GA runs is depicted in Figure~\ref{fig:pareto_front}. See Appendix \ref{app:testfunc} for how we validated the maximum generation number. The convergence of the front is confirmed by the large number of generations with the population size of 372, involving around 1.1 million XFOIL evaluations of $CLD_{max}$ and $\Delta \alpha$. Without duplicates, a set of 80 Pareto-optimal airfoil shapes was obtained via DbM from the 25 chosen baseline shapes. For comparison, the whole UIUC database \parencite{uiuc_database}, as well as the baseline cases, are evaluated and plotted in Figure~\ref{fig:pareto_front} together. It is noted that baseline \#19 has zero $CLD_{max}$ and $\Delta \alpha$ because it is intentionally mirrored, and XFOIL failed to evaluate its aerodynamic performance. We assigned zero values to failing cases like this because they represented airfoil geometries found to be aerodynamically unviable in the XFOIL space. The GA optimization successfully developed the Pareto-front, with two ends at $\left(CLD_{max}, \Delta \alpha \right) = \left( 30.63, 40^{\circ}\right)$ and $\left(CLD_{max}, \Delta \alpha \right) = \left( 264.17, 11^{\circ}\right)$. Even in the largest maximum lift-drag ratio case, the angle of attack gap between the stall and design point is $11^{\circ}$, providing the airfoil with a decent tolerant range for off-design operations.

The Pareto-front is divided into three different clusters, each constitutes a segment of the front that does not overlap with the others. It is worth noting that the non-overlapping division of the front is a result of clustering through Principal Component Analysis (PCA), rather than being manually assigned. The details of the clustering are provided in Appendix~\ref{app:pca}.

Figure~\ref{fig:pareto_shape} depicts nine representative optimal airfoil shapes on the Pareto-front, arranged in ascending order of $CLD_{max}$. From each cluster, three airfoil shapes with distinct objective function values have been selected for representation. Also, note that Figure~\ref{fig:pareto_shape}\subref{fig:pareto_shape_a} shows the extreme case of the smallest $CLD_{max}$ and largest $\Delta \alpha$, while Figure~\ref{fig:pareto_shape}\subref{fig:pareto_shape_i} depicts the opposite extreme of the largest $CLD_{max}$ and smallest $\Delta \alpha$. It can be seen that within each cluster the overall shape remains unchanged, with only a gradual decrease in airfoil thickness as $CLD_{max}$ increases. Since thin airfoils such as bird-like airfoils \parencite{ananda2018design}, e.g. \#13 and \#14 of the baseline shapes, are known for their high $CLD$ performance, the trend of airfoil thickness observed in the Pareto-front appears to be reasonable.

\begin{figure*}\centering
\begin{subfigure}[t]{0.33\linewidth}
\vbox{
\centering{
  \includegraphics[width=\linewidth]{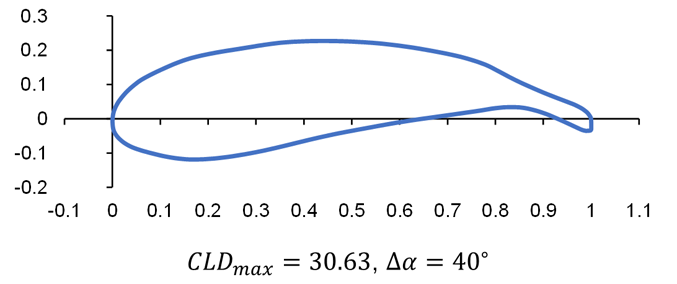}
}%
}%
\subcaption{\centering\label{fig:pareto_shape_a}}
\end{subfigure}
\begin{subfigure}[t]{0.33\linewidth}
\centering{
\includegraphics[width=\linewidth]{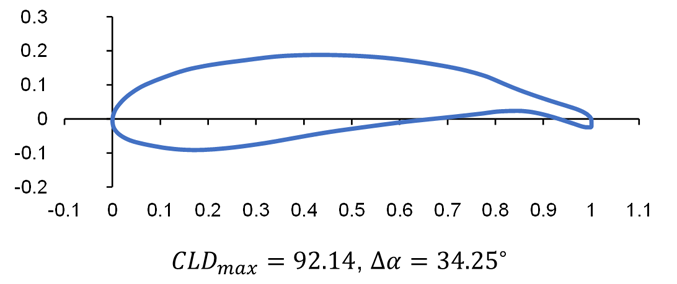}
\subcaption{\centering\label{fig:pareto_shape_b}}
}\end{subfigure}%
\begin{subfigure}[t]{0.33\linewidth}
\vbox{
\centering{
\includegraphics[width=\linewidth]{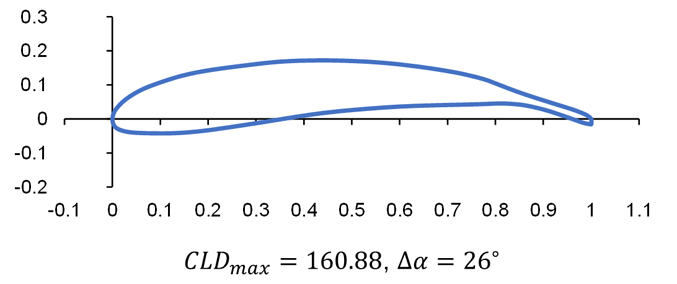}
}%
}
\subcaption{\centering\label{fig:pareto_shape_c}}
\end{subfigure}\\
\begin{subfigure}[t]{0.33\linewidth}
\vbox{
\vspace*{1.7em}%
\centering{
  \includegraphics[width=\linewidth]{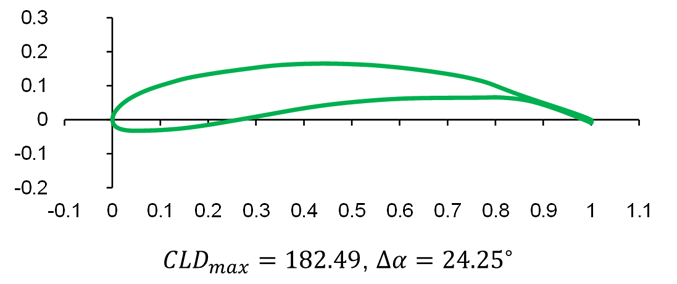}
}%
}%
\subcaption{\centering\label{fig:pareto_shape_d}}
\end{subfigure}
\begin{subfigure}[t]{0.33\linewidth}
\centering{
\includegraphics[width=\linewidth]{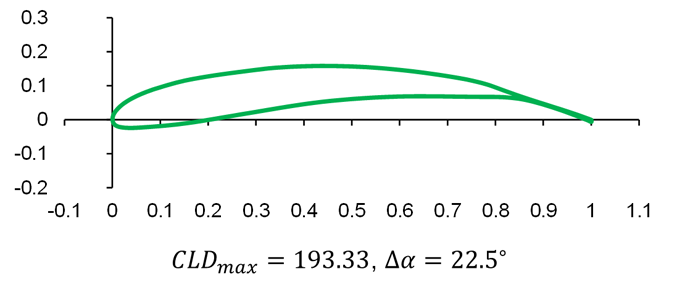}
\subcaption{\centering\label{fig:pareto_shape_e}}
}\end{subfigure}%
\begin{subfigure}[t]{0.33\linewidth}
\vbox{
\vspace*{1.7em}%
\centering{
\includegraphics[width=\linewidth]{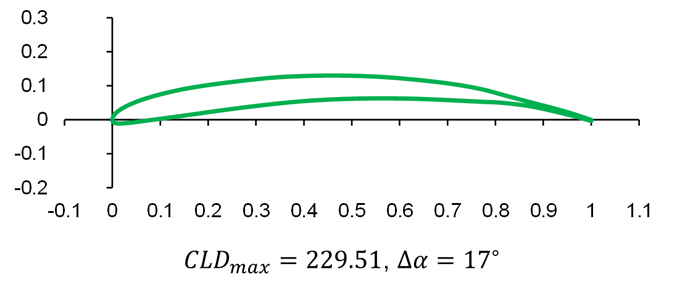}
}%
}
\subcaption{\centering\label{fig:pareto_shape_f}}
\end{subfigure}\\
\begin{subfigure}[t]{0.33\linewidth}
\vbox{
\vspace*{1.7em}%
\centering{
  \includegraphics[width=\linewidth]{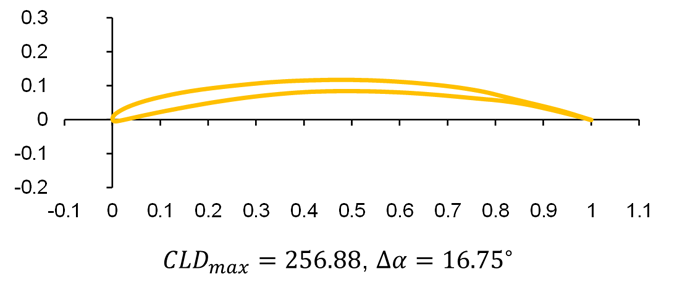}
}%
}%
\subcaption{\centering\label{fig:pareto_shape_g}}
\end{subfigure}
\begin{subfigure}[t]{0.33\linewidth}
\centering{
\includegraphics[width=\linewidth]{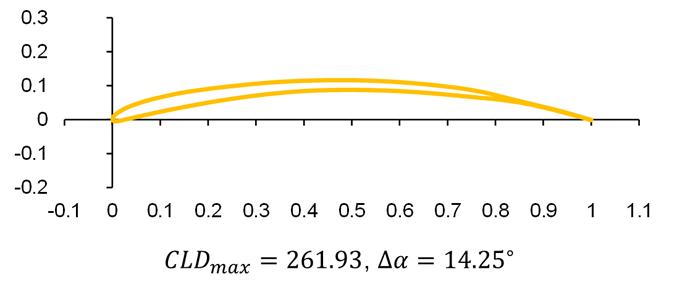}
\subcaption{\centering\label{fig:pareto_shape_h}}
}\end{subfigure}%
\begin{subfigure}[t]{0.33\linewidth}
\vbox{
\vspace*{1.7em}%
\centering{
\includegraphics[width=\linewidth]{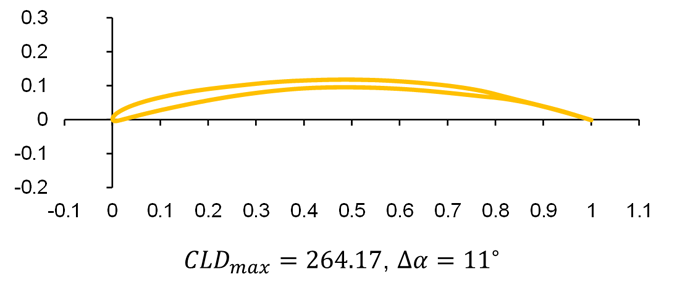}
}%
}
\subcaption{\centering\label{fig:pareto_shape_i}}
\end{subfigure}
\caption{Nine representative Pareto-optimal airfoil shapes. (\subref{fig:pareto_shape_a})-(\subref{fig:pareto_shape_c}) are in cluster 1, (\subref{fig:pareto_shape_d})-(\subref{fig:pareto_shape_f}) are in cluster 2 and (\subref{fig:pareto_shape_g})-(\subref{fig:pareto_shape_i}) are in cluster 3. See Appendix \ref{app:pca} for the clustering.}\label{fig:pareto_shape}
\end{figure*}

\begin{figure}[tb!]\centering
\begin{subfigure}[t]{\linewidth}
\vbox{
\centering{
  \includegraphics[width=.85\linewidth]{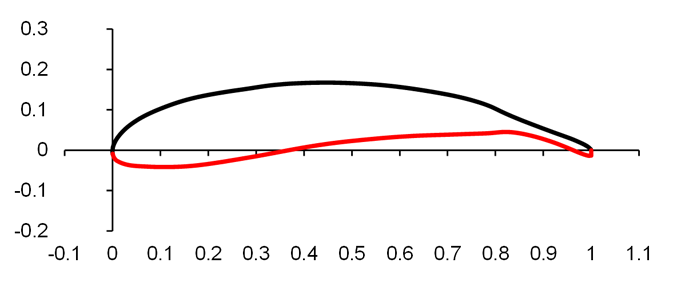}
}%
}%
\subcaption{\centering\label{fig:pca_shape_a} Total mean}
\end{subfigure}\\
\begin{subfigure}[t]{\linewidth}
\vbox{
\centering{
  \includegraphics[width=.85\linewidth]{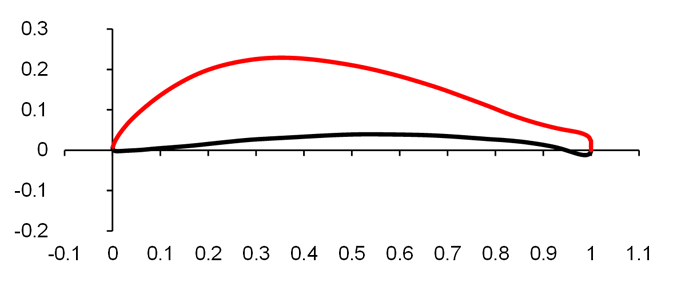}
}%
}%
\subcaption{\centering\label{fig:pca_shape_b} PCA axis 1}
\end{subfigure}\\
\begin{subfigure}[t]{\linewidth}
\vbox{
\centering{
  \includegraphics[width=.85\linewidth]{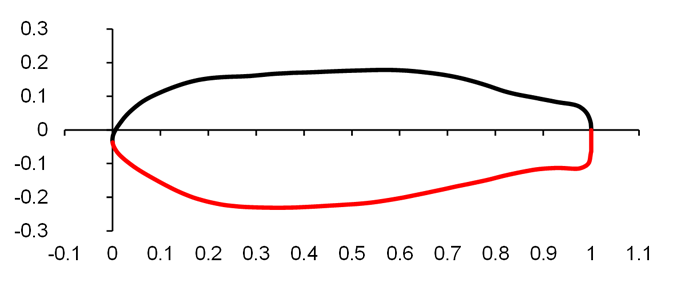}
}%
}%
\subcaption{\centering\label{fig:pca_shape_c} PCA axis 2}
\end{subfigure}
\caption{Morphed airfoil shapes generated by the optimal weight vectors, representing (\subref{fig:pca_shape_a}) the total mean of all optimal airfoils' weights, (\subref{fig:pca_shape_b}) the coefficients of the PCA axis having the most variance and (\subref{fig:pca_shape_c}) the coefficients of the PCA axis having the second-most variance. The black and red surfaces correspond to the first and second half of the collocation points, respectively.}\label{fig:pca_shape}
\end{figure}

Cluster 1, made up by 48 optimal airfoil shapes, resembles the total mean of the Pareto-front, which is the average of all airfoil shapes on the Pareto-front (see Figure~\ref{fig:pca_shape}\subref{fig:pca_shape_a}). This makes sense as they make up the majority of airfoil shapes located on the front. Moreover, this cluster is located near the origin in the PCA-projected weight space (see Figure~\ref{fig:pca_cluster} in Appendix \ref{app:pca}), indicating that no radical morphing of the airfoil shape took place from the mean shape.

Next, cluster 2 contains 24 optimal airfoil shapes. Compared to those in cluster 1, the most distinguishing feature is their narrow trailing edge regions, which are typically favorable for increasing lift. However, these airfoils are not greatly different from the origin in the PCA-projected weight space and are close to the total mean Pareto-front.

Finally, 8 optimal airfoil shapes are found in cluster 3 from the optimization. This cluster includes the airfoil shapes experiencing more drastic morphing than the other clusters. This is manifested by the fact that they are the thinnest airfoils where the leading edge region's thickness also becomes narrow.

\begin{figure*}
  \centering
  \includegraphics[width=\linewidth]{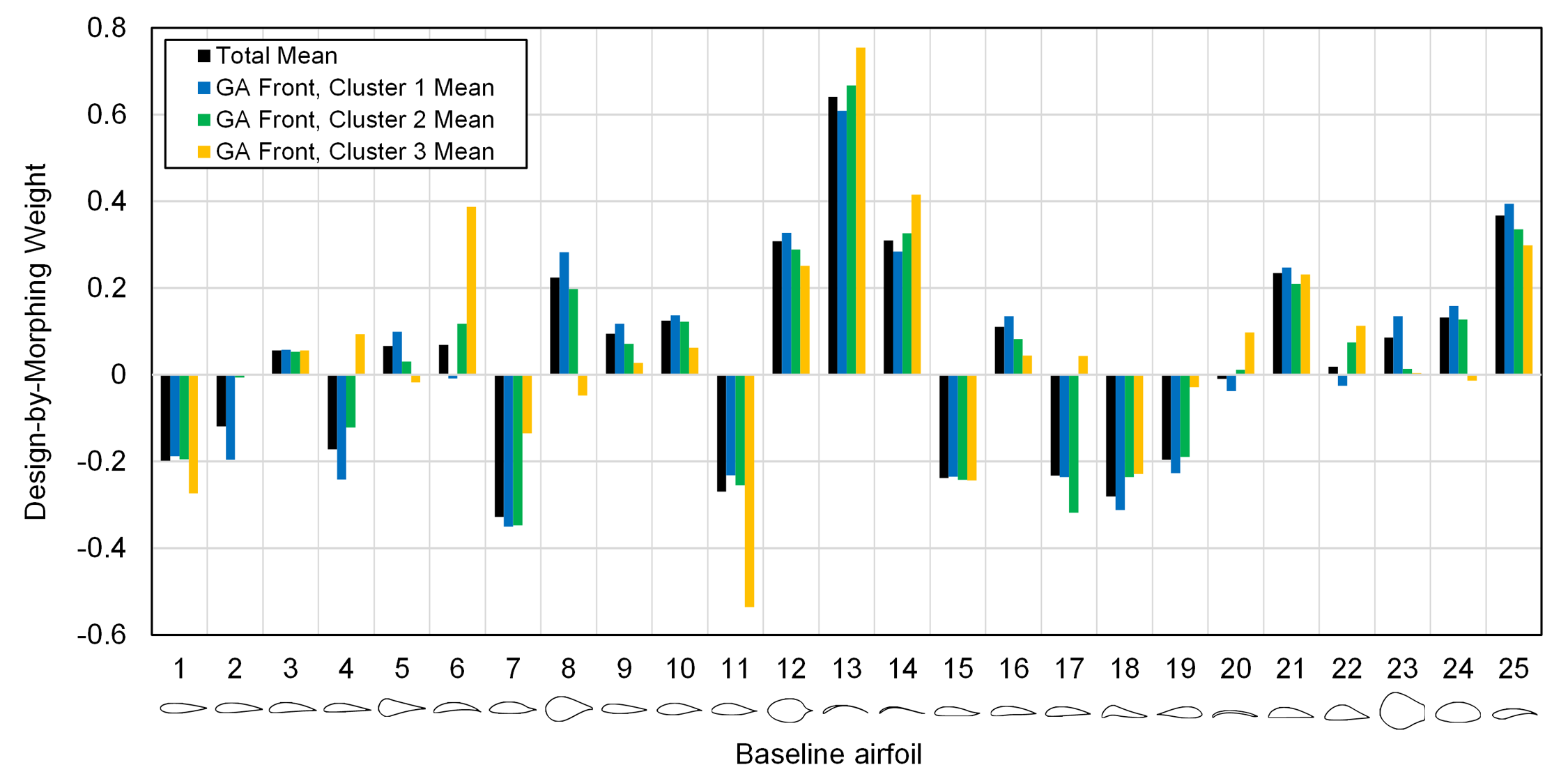}
  \caption{Mean weight distributions of the Pareto-optimal airfoil shapes with respect to twenty-five baseline airfoil shapes.}
  \label{fig:weightdist}
\end{figure*}

The mean weight distributions with respect to 25 original baseline shapes are shown in Figure~\ref{fig:weightdist}. Overall, the weight distributions of the three clusters conform to the weight distribution of the total mean. It turned out that baseline shape \#13 (model name: AS6097) was commonly the most significant for morphing. Since this baseline shape has the best in $CLD_{max}$ and the second best in $\Delta \alpha$ among the 25 baseline shapes (see Figure~\ref{fig:pareto_front}), it was likely to persist in the GA runs over generations against the selection pressure that only sorts out dominant individuals in terms of both $CLD_{max}$ and $\Delta \alpha$. However, excellence in the objectives of an individual baseline shape does not necessarily guarantee its survival, which is the case for the globally best baseline shape \#6 (model name: AH 79-100C). An individual's superior `phenotype' may be no longer revealed or even suppressed after the morphing is done and all `genes' are mixed up. In the same sense, inferiority in the objectives of an individual does not necessarily result in elimination, as demonstrated by the `mirrored' baseline \#19.

As discussed from the examination of the morphed airfoil shapes, both cluster 1 and 2 show no significant differences from the total mean Pareto-front. Through small shape variations in shape from the total mean Pareto-front, as seen in Figure~\ref{fig:pca_shape}\subref{fig:pca_shape_a}, reaching these optima would be relatively easily. In contrast, cluster 3 has a number of weights that are quite different from the mean (e.g., \#6 and \#11) and substantial morphing would be required if one starts with the total mean airfoil shape.

In the context of the present study, each axis obtained by PCA can be considered as a unique morphed airfoil shape because the 25 PCA coefficient vectors defined in the weight space $\mathcal{D}^{25}$ are orthogonal to each other. These 25 new morphed airfoils span the entire design space and therefore can serve as alternative baseline shapes in lieu of the original ones. More importantly, the dominance of the first two PCA axes with respect to the data point variance, accounting for 95\% of total variance explained, suggests that the major geometric features of the 208 airfoil shapes we found through optimization are virtually generated by morphing of these two new airfoils. A small variance of a PCA axis indicates that the data points are not significantly deviated from their mean on the axis. In other words, the baseline shape corresponding to this PCA axis has an marginal impact on morphing the airfoil shape for optimization.

Once we pick two baseline shapes from the first two dominant PCA axes, whose associated collocation vectors are denoted as $\vec{P}_1$ and $\vec{P}_2$ for example, and use them to morph the airfoil shape obtained from the total mean of the Pareto-optimal weight vector set, which corresponds to the mean collocation vector $\vec{P}_{mean}$, we gain a better understanding of how the morphing, especially along each PCA axis, influences major geometric changes in the optimal airfoil shapes. These airfoil shapes are depicted in Figure~\ref{fig:pca_shape}, where the black and red surfaces are distinguished to emphasize that they represent the first and second halves of the collocation points, respectively. For example, we note that the orientation of the two surface of $\vec{P}_1$ is mirrored in comparison to that of $\vec{P}_{mean}$, meaning that the stronger the weight of PCA axis 1 in the positive direction is, the narrower the morphed airfoil shape gets.

\section{Discussion}

Most shape parameterization methods rely upon the careful selection of geometric constraints and parameters, which determines the likelihood of success. The fidelity offered by such methods largely depends on the number of the parameters chosen. Moreover, these designs are limited by the parametric constraints and the implicit designer bias, making it difficult to perform extrapolation or make radical global changes. Data-driven methods typically rely on the assumption that the optimum solutions are not far from the training data-set, which again limits the ability to make radical shape changes.

DbM, on the other hand, creates a design space that is not inhibited by geometric constraints, enables extrapolation from the design space, which is particularly useful for airfoil design, and is applicable to a wide range of engineering design problems. It does not suffer from the curse of dimensionality when parameterizing airfoils by control points and allows for high-fidelity representation of airfoils without increasing the number of independent parameters in the problem. Using only 25 baseline shapes from the UIUC database, we were able to recreate 99.87\% of UIUC database with an MAE error < 1\%. We also showed that extraordinary and broad searches are possible using DbM. By applying it to the bi-objective shape optimization with the objectives of maximizing $CLD_{max}$ and $\Delta \alpha$, we could achieve significant results compared to our baseline shapes. We posit that for the design parameterization of airfoils as well as other 2D/3D shapes, DbM should be the preferred method for creating an unconstrained, unbiased and non-database-driven design space that allows for radical modifications, which can often result in non-conforming shapes.
In this paper, our qualitative selection of 25 baseline shapes adequately spanned the design space with tolerable error. However, it would be possible that even a smaller number of baseline shapes than 25 could successfully construct the design space if some of the current baseline shapes were redundant. To further understand DbM, our future work will focus on performing sensitivity analysis of DbM on the baseline shape selection and applying DbM for design of turbo-machinery.

\section{Conclusion}
The DbM design strategy creates a design space that contains novel and radical 2D airfoils that are not constrained by geometric parameters or designer bias. Optimization within the design space created, for the dual objectives of $CLD_{max}$ and $\Delta \alpha$, resulted in remarkable improvements in both objectives and provided a Pareto-front of optimal airfoil designs. The final airfoils showed significant advancements compared to the input baseline shapes. 

Overall, in our optimization study with respect to the 2D airfoil optimization problem, DbM is a suitable method for design space creation. In addition, our methodology is highly adaptable and can be utilized for shape optimization of other fluid machinery. Our ongoing work includes the applications of DbM in conjunction with Bayesian optimization to 3D airfoil optimization and vertical-axis wind turbine optimization problems.
\section*{Acknowledgment} 

The authors would like to thank Prof. \"{O}mer Savas\c{s}, affiliated with University of California at Berkeley, for providing insightful discussions regarding airfoils and aerodynamics. The authors acknowledge the use of the Extreme Science and Engineering Discovery Environment (XSEDE), supported by National Science Foundation grant number ACI-1548562 through allocation TG-CTS190047.


\section*{Data and Materials Availability} 
\noindent
The data needed to evaluate the conclusions are present in the paper and Appendices. The data files and optimization setup will be posted in a public repository upon publication of the paper.
\nomenclature[A]{$c$}{Airfoil chord length ($\mathrm{m}$)}
\nomenclature[A]{$d$}{Drag force of an airfoil per unit span ($\mathrm{N~m^{-1}}$)}
\nomenclature[A]{$l$}{Lift force of an airfoil per unit span ($\mathrm{N~m^{-1}}$)}
\nomenclature[A]{$\vec{P}$}{$y$-coordinate collocation vector of a morphed airfoil}
\nomenclature[A]{$\vec{S}$}{$y$-coordinate collocation vector of a baseline airfoil}
\nomenclature[A]{$U$}{Free-stream flow speed ($\mathrm{m~s^{-1}}$)}
\nomenclature[A]{$w$}{Design-by-Morphing weight factor}

\nomenclature[B]{$\alpha$}{Airfoil angle of attack ($\mathrm{^{\circ}}$)}
\nomenclature[B]{$\alpha_s$}{Airfoil stall angle ($\mathrm{^{\circ}}$)}
\nomenclature[B]{$\Delta \alpha$}{Stall angle tolerance, the range of $\alpha$ between the stall point and the maximum lift-drag ratio point ($\mathrm{^{\circ}}$)}
\nomenclature[B]{$\nu$}{Fluid kinematic viscosity ($\mathrm{m^2~s^{-1}}$)}
\nomenclature[B]{$\rho$}{Fluid density ($\mathrm{kg~m^{-3}}$)}

\nomenclature[C]{$C_d$}{Drag coefficient of an airfoil per unit span, $2l/(\rho U^2 c)$}
\nomenclature[C]{$C_l$}{Lift coefficient of an airfoil per unit span, $2d/(\rho U^2 c)$}
\nomenclature[C]{$CLD$}{Lift-drag ratio of an airfoil, $C_l / C_d$}
\nomenclature[C]{$CLD_{max}$}{Maximum lift-drag ratio of an airfoil, $\max_{\alpha} CLD (\alpha)$}
\nomenclature[C]{$\mathrm{Re}$}{Reynolds number based on airfoil chord length, $U c/\nu$}
\printnomenclature

\appendix   
\section{Aerodynamic Optimization Objectives}\label{app:prelim}

Airfoil optimization has become a common practice in aerodynamic design problems that involve maximization of one or more performance parameters of airfoils. We mainly consider the following two performance parameters: the lift-drag ratio and stall angle. Given the flow speed $U$, fluid density $\rho$ and airfoil chord length $c$, the lift and drag coefficients of an airfoil per unit span at an angle of attack $\alpha$, $C_l$ and $C_d$, are expressed as:
\begin{equation}
C_l(\alpha) \equiv \frac{l(\alpha)}{\frac{1}{2}\rho U^2 c}
~,~~~
C_d(\alpha) \equiv \frac{d(\alpha)}{\frac{1}{2}\rho U^2 c}.
\label{eq:clcd}
\end{equation}
where $l$ and $d$ are lift and drag force per unit span, respectively, both of which change with respect to $\alpha$. In this paper, these parameters are predicted using XFOIL \parencite{drela_1989}, a program for analyzing a subsonic 2D airfoil, with varying $\alpha$ and then used for optimization. Based on $C_l$ and $C_d$, the lift-drag ratio $CLD$ is calculated as:
\begin{equation}
CLD(\alpha) = \frac{C_l \,(\alpha)}{C_d \,(\alpha)}.
\label{eq:l/d}
\end{equation}
On the other hand, we define the stall angle $\alpha_s$ as an angle of attack where $C_l$ reaches its first local maximum as the angle increases from 0${}^{\circ}$, or:
\begin{equation}
\begin{gathered}
\alpha_s \equiv \min_{\alpha \ge 0} \alpha ~~ \text{where}~~ {}^{\exists} \delta > 0 ~~ \text{such that}\\ C_l\,(\alpha) \ge C_l\,(x) ~~ {}^{\forall} x \in [\alpha - \delta ,~ \alpha + \delta] .
\end{gathered}
\label{eq:stall}
\end{equation}
Note that this definition is more conservative than the typical definition of stalling, where the flow at the rear region begins to fully separate and $C_l$ is globally maximized. $\alpha_s$ is occasionally smaller than the global maximum of $C_l$. Nonetheless, this approach helps avoid overestimation of the stall angle, which is expected to occur in XFOIL due to the nature of its flow solver having a limited accuracy in stall and post-stall conditions.

\begin{figure}[tb!]
    \vspace*{-1.5em}%
    \centering
    \includegraphics[width=\linewidth]{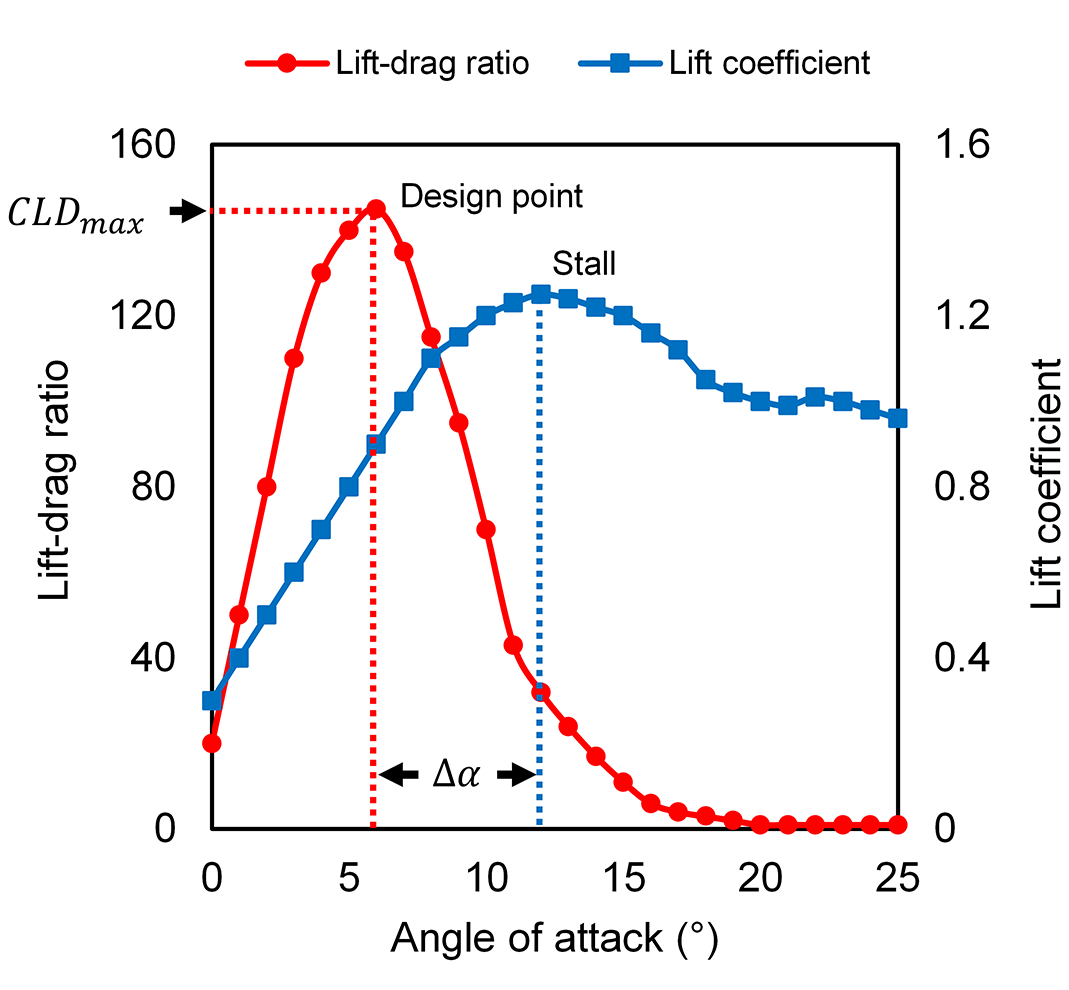}
    \caption{Schematic diagram of airfoil performance curves.} 
    \label{fig:prelim}
\end{figure}

$CLD$ and $\alpha_s$ have been typically considered to be significant in characterizing airfoil performance. For example, when it comes to lift-type wind turbines, the point where $CLD$ is maximized is often chosen as the design point. However, since wind turbines cannot always operate under design conditions, $\alpha_s$ needs to be additionally considered to evaluate how far they run under increasing lift conditions. For well-designed airfoils, $\alpha_s$ generally occurs later than the design point, which yields operational tolerance beyond the design point. Consequently, the stall angle tolerance $\Delta \alpha$, i.e. the range between these two angles of attack, expressed as
\begin{equation}
\Delta \alpha \equiv \max \left(0,~ \alpha_s - \argmax_{\alpha \in \mathbb{R}}{CLD (\alpha) } \right) ,
\label{eq:deltaalpha}
\end{equation}
can be a proper choice to evaluate the off-design performance  \parencite{li_method_2013}. Figure~\ref{fig:prelim} depicts a schematic diagram of how $CLD$ and $\Delta \alpha$ are determined on airfoil performance curves.

\section{Baseline Airfoil Shapes and Evaluation} \label{app:baseline}

\begin{table*}[hbt!]
\caption{The model names, features, shape outlines, and XFOIL evaluation results of the 25 baseline shapes used by DbM in this paper. The coordinates of the baseline shapes are obtained from the UIUC airfoil coordinates database \parencite{uiuc_database}. The airfoil evaluation results are obtained for an incompressible flow with a chord Reynold number of $1\times10^6$. The reference evaluation results are interpolated from the online XFOIL database \parencite{airfoiltools}; N/A indicates that no data is available.}
\label{tab:baseshapetab}%
\centering{%
\begin{tabularx}{\linewidth}{c C C m{0.7in} c c c c}
\toprule
\multirow{2}{*}{Index} &
  \multirow{2}{*}{Model Name} &
  \multirow{2}{*}{Series (Features)} &
  \multicolumn{1}{c}{\multirow{2}{*}{Airfoil Shape}} &
  \multicolumn{2}{c}{\small Reference} &
  \multicolumn{2}{c}{\small Present} \\
   &                        &                         &  & $CLD_{max}$ & $\Delta\alpha$ & $CLD_{max}$ & $\Delta\alpha$ \\ \midrule

1  & NACA 0012              & NACA (4-digit)          &\includegraphics[scale=0.5]{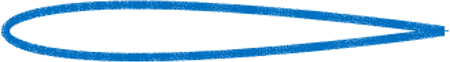}  & 75.6     & 8.50   & 69.3     & 6.75   \\
2  & NACA 2412              & NACA (4-digit)          &\includegraphics[scale=0.5]{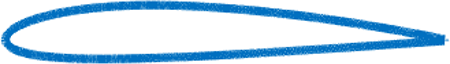}  & 101.4    & 12.00  & 99.5     & 12.00  \\
3  & NACA 4412              & NACA (4-digit)          &\includegraphics[scale=0.5]{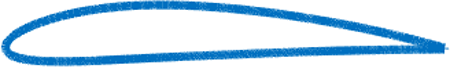}  & 129.4    & 1.75   & 126.2    & 11.50  \\
4  & E 205                  & Eppler                  &\includegraphics[scale=0.5]{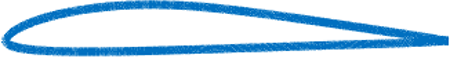}  & 128.3    & 8.50   & 124.4    & 10.50  \\
5  & AH 81-K-144 W-F Klappe & Althaus                 &\includegraphics[scale=0.5]{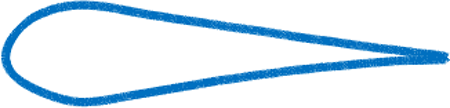}  & 89.7     & 2.00   & 91.6     & 2.00   \\
6  & AH 79-100 C            & Althaus                 &\includegraphics[scale=0.5]{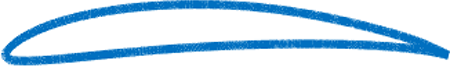}  & 183.0    & 14.75  & 170.6    & 15.50  \\
7  & AH 79-K-143/18         & Althaus                 &\includegraphics[scale=0.5]{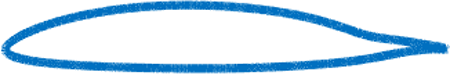}  & 110.9    & 1.50   & 107.0    & 1.50   \\
8  & AH 94-W-301            & Althaus                 &\includegraphics[scale=0.5]{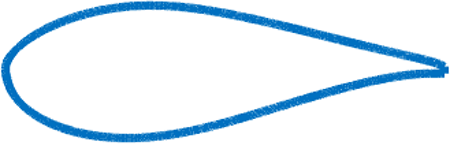}  & 103.0    & 4.00   & 101.4    & 2.75   \\
9  & NACA 23112             & NACA (5-digit)          &\includegraphics[scale=0.5]{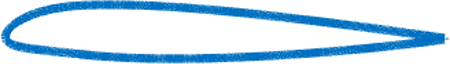}  & 98.6     & 6.75   & 96.9     & 8.00   \\
10 & NACA 64(2)-415         & NACA (6-digit)          &\includegraphics[scale=0.5]{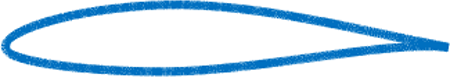}  & 120.6    & 12.50  & 113.8    & 13.00  \\
11 & NACA 747(A)-315        & NACA (7-digit)          &\includegraphics[scale=0.5]{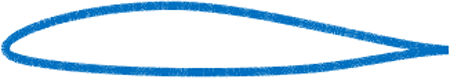}  & 111.5    & 12.00  & 105.8    & 13.00  \\
12 & Griffith 30\% Suction  & Griffith (Suction)      &\includegraphics[scale=0.5]{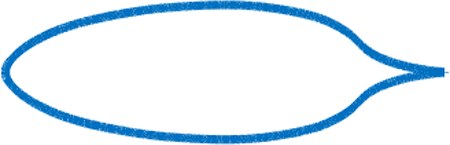}  & 17.3     & 0.00   & 17.9     & 0.00   \\
13 & AS 6097                & Selig (Bird-like)       &\includegraphics[scale=0.5]{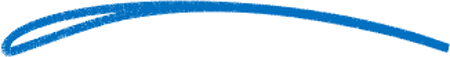}  & N/A      & N/A    & 171.2    & 14.00  \\
14 & E 379                  & Eppler (Bird-like)      &\includegraphics[scale=0.5]{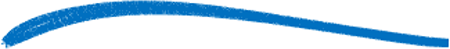}  & N/A      & N/A    & 160.0    & 2.00   \\
15 & Clark YS               & Clark                   &\includegraphics[scale=0.5]{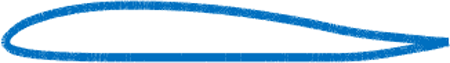}  & 85.7     & 5.25   & 82.3     & 5.75   \\
16 & Clark W                & Clark                   &\includegraphics[scale=0.5]{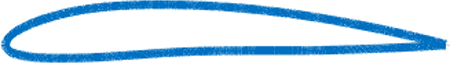}  & 116.1    & 11.00  & 114.8    & 11.00  \\
17 & Clark Y                & Clark                   &\includegraphics[scale=0.5]{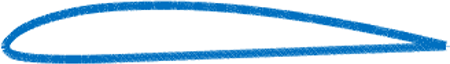}  & 114.8    & 11.75  & 113.7    & 12.75  \\
18 & Chen                   & Chen                    &\includegraphics[scale=0.5]{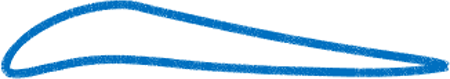}  & 125.4    & 0.00   & 126.7    & 0.00   \\
19 & S2027 Mirrored          & Selig (Mirrored)         &\includegraphics[scale=0.5]{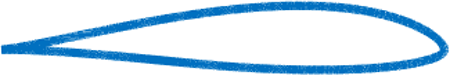}  & N/A      & N/A    & 0.00      & 0.00    \\
20 & GOE 417A               & Gottingen (Thin plate)  &\includegraphics[scale=0.5]{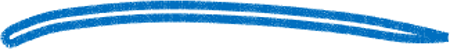}  & 86.7     & 5.25   & 90.4     & 5.25   \\
21 & GOE 611                & Gottingen (Flat bottom) &\includegraphics[scale=0.5]{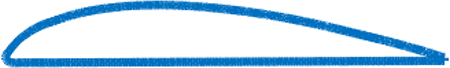}  & 125.6    & 9.00   & 129.7    & 9.00   \\
22 & Dragonfly Canard       & Dragonfly               &\includegraphics[scale=0.5]{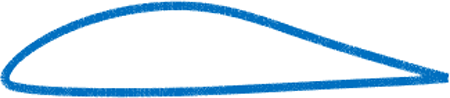}  & 144.6    & 2.50   & 147.5    & 3.00   \\
23 & FX 79-W-470A           & Wortmann (Fat)          &\includegraphics[scale=0.5]{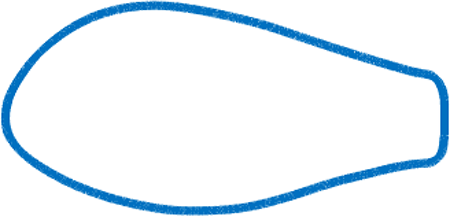}  & N/A      & N/A    & 23.9     & 9.25   \\
24 & Sikorsky DBLN-526      & Sikorsky (Fat)          &\includegraphics[scale=0.5]{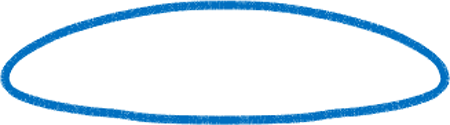}  & 53.3     & 4.75   & 51.5     & 4.25   \\
25 & FX 82-512              & Wortmann                &\includegraphics[scale=0.5]{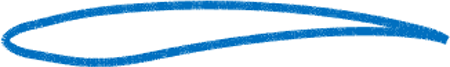}  & 99.1     & 14.75  & 98.7     & 13.00  \\ \bottomrule
\end{tabularx}
}%
\end{table*}

Our optimization methodology does not rely on one specific airfoil evaluation tool. To compare our results with previous literature and help future researchers quickly reproduce our results, we used XFOIL \parencite{drela_1989} in the present study. The two design objectives, $CLD_{max}$ and $\Delta\alpha$, are obtained from the $C_l$ and $C_d$ data calculated by the XFOIL at different angle of attacks (see Figure~\ref{fig:prelim}). 

For improved efficiency and consistency, we used XFOIL to generate performance data and did not rely on any of its built-in paneling features. The conditioning and re-paneling of the morphed airfoil coordinates are custom-built at the end of our DbM algorithm, transforming the coordinates into 200 or 250 vortex panels with a relatively higher concentration where the curvature is high. To reduce evaluation time, we first performed a rough scan with an $\alpha$ increment of $1^{\circ}$ to estimate the range determining $\Delta \alpha$, and then finer scans for $CLD_{max}$ and $\Delta \alpha$ separately with an $\alpha$ increment of $0.25^{\circ}$ within and around the estimated range of $\Delta \alpha$ from the initial rough scan.

It is worth noting that XFOIL uses a global Newton's method \parencite{drela_1989} to solve the boundary layer and transition equations simultaneously and uses the solution from the previous angle of attack as a starting guess. As a result, ill-conditioned airfoil coordinates and the occurrence of flow separation can both lead to non-convergence of the XFOIL evaluation. To ensure the robustness and correctness of our airfoil evaluation, our XFOIL wrapper attempts to reach convergence by restarting the root-finding with a fresh starting guess and gradually increasing the number of panels. If both attempts fail, the wrapper will check convergence at neighboring points, which will indicate whether flow separation occurs or not. Besides non-convergence, we further verify the correctness of the Newton's method by comparing the calculated viscous and inviscid drag coefficients. The later is determined purely by the potential flow theorem and have to be smaller than its viscous counterpart due to its neglect of the friction (viscous effect). Any angle with incorrect result will undergo the same treatment as non-converging ones, hence ensuring the correctness of our airfoil performance evaluation. A comparison between our XFOIL evaluation and an existing database of the same airfoils under the same flow conditions is provided in Table \ref{tab:baseshapetab}.

\section{Optimization Test Functions and Validation} \label{app:testfunc}

\begin{figure}[tb!]\centering
\begin{subfigure}[tb!]{\linewidth}
\centering{\hspace*{.3in}
  \includegraphics[width=1.6in]{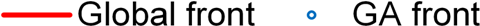}
}%
\end{subfigure}\\[-.5em]
\begin{subfigure}[t]{.49\linewidth}\captionsetup{justification=raggedleft,singlelinecheck=false}
\vbox{
\vspace*{1.7em}%
\centering{
  \includegraphics[width=\linewidth]{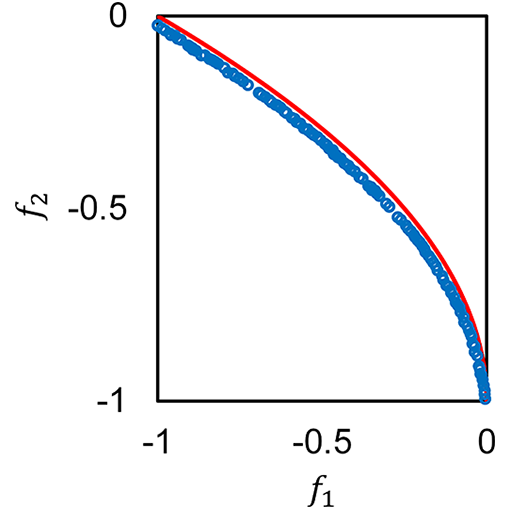}
}%
}%
\subcaption{\label{fig:benchmark_a} ZDT1 \hspace*{.4in}}
\end{subfigure}
\begin{subfigure}[t]{.49\linewidth}\captionsetup{justification=raggedleft,singlelinecheck=false}
\centering{
\includegraphics[width=\linewidth]{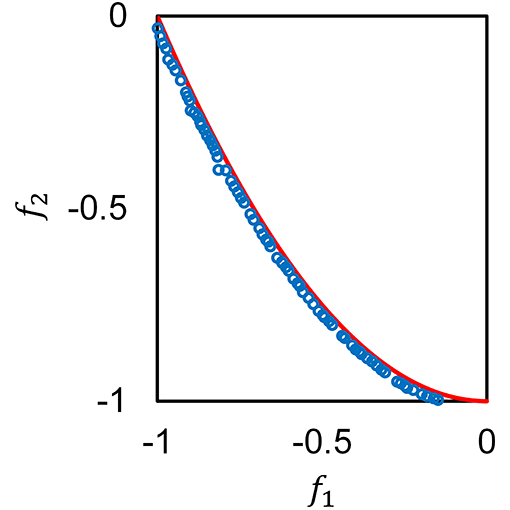}
\subcaption{\label{fig:benchmark_b} ZDT2 \hspace*{.4in}}
}\end{subfigure}\\
\begin{subfigure}[t]{.49\linewidth}\captionsetup{justification=raggedleft,singlelinecheck=false}
\vbox{
\vspace*{1.7em}%
\centering{
  \includegraphics[width=\linewidth]{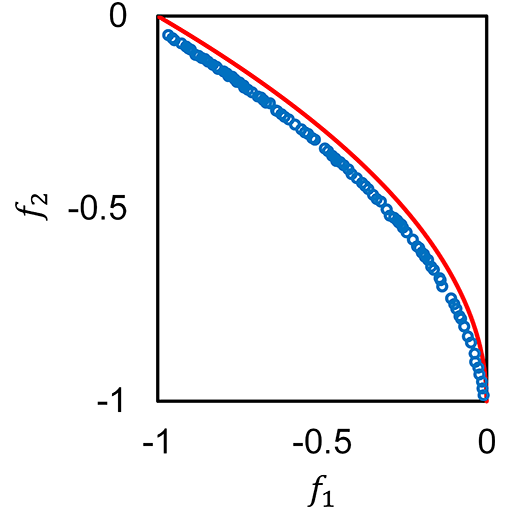}
}%
}%
\subcaption{\label{fig:benchmark_c} ZDT4 \hspace*{.4in}}
\end{subfigure}
\begin{subfigure}[t]{.49\linewidth}\captionsetup{justification=raggedleft,singlelinecheck=false}
\centering{  \hspace*{0.025in}
\includegraphics[width=\linewidth]{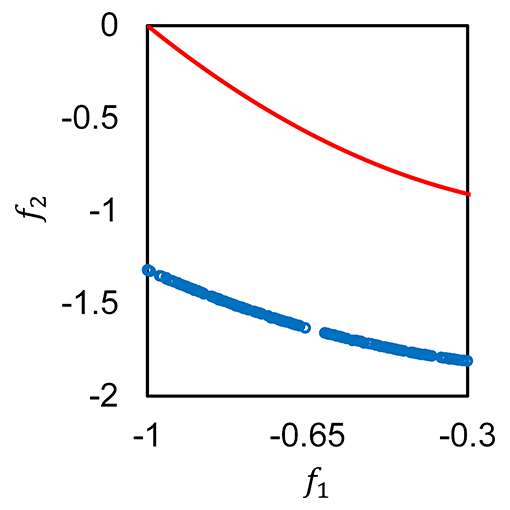}
\subcaption{\label{fig:benchmark_d} ZDT6 \hspace*{.4in}}
}\end{subfigure}
\caption{Multi-objective optimization of benchmark test functions using GA.}\label{fig:benchmark}
\end{figure}


We used the multi-objective problems suggested by \textcite{10.1162/106365600568202} to test our GA setup. The details of the test functions are given in Table \ref{table:a2}. All the test functions are aimed to be minimized with 25 variables in the given design space.

\begin{table*}[]
\caption{Benchmark Test Functions. All of the test functions are bi-objective with extended to \textit{n}-dimensional constrained search space.}\label{table:a2}
\centering{%
  \begin{tabularx}{\linewidth}{c C C C c}
    \toprule 
   Problem & Bounds & Objective Functions & Optima & Note\\\midrule
    $ZDT1$ & $\begin{aligned} w_i & \in [0,1], \\[-.5ex] i &= 1,\dots,n \end{aligned}$   & $\begin{aligned}  f_1(\vec{w})&=w_1 \\[-.5ex] f_2(\vec{w})&=g(\vec{w})\left[1-(f_1(\vec{w})/g(\vec{w}))^{1/2}\right]\\[-.5ex] g(\vec{w})&=1+9\left(\sum_{i=2}^nw_i\right)/(n-1) \end{aligned}$   & $\begin{aligned} w_1 & \in [0,1]\\[-.5ex] w_i &=0, \\[-.5ex] i &= 2,\dots,n \end{aligned}$  & convex\\ \\

    $ZDT2$ & $\begin{aligned}  w_i & \in [0,1], \\[-.5ex] i &= 1,\dots,n \end{aligned}$  & $\begin{aligned}  f_1(\vec{w})&=w_1 \\[-.5ex] f_2(\vec{w})&=g(\vec{w})\left[1-(f_1(\vec{w})/g(\vec{w}))^{2}\right]\\[-.5ex] g(\vec{w})&=1+9\left(\sum_{i=2}^nw_i\right)/(n-1) \end{aligned}$   & $\begin{aligned} w_1 & \in [0,1]\\[-.5ex] w_i &=0, \\[-.5ex] i &= 2,\dots,n \end{aligned}$  & non-convex\\ \\

    $ZDT4$ & $\begin{aligned} w_1 & \in [0,1]\\[-.5ex] w_i & \in [-5,5], \\[-.5ex] i &= 2,\dots,n \end{aligned}$   & $\begin{aligned}  f_1(\vec{w})&=w_1 \\[-.5ex] f_2(\vec{w})&=g(\vec{w})\left[1-(f_1(\vec{w})/g(\vec{w}))^{1/2}\right]\\[-.5ex] g(\vec{w})&=10n+\sum_{i=2}^n\left(w_i^2-10\cos(4\pi w_i)\right)-9 \end{aligned}$   & $\begin{aligned} w_1 & \in [0,1]\\[-.5ex] w_i &=0,  \\[-.5ex] i &= 2,\dots,n \end{aligned}$  & non-convex\\ \\

    $ZDT6$ & $\begin{aligned}  w_i & \in [0,1], \\[-.5ex] i &= 1,\dots,n \end{aligned}$  & $\begin{aligned}  f_1(\vec{w})&=1-\exp(-4w_1)\sin^6(6 \pi w_1) \\[-.5ex] f_2(\vec{w})&=g(\vec{w})\left[1-(f_1(\vec{w})/g(\vec{w}))^{2}\right]\\[-.5ex] g(\vec{w})&=1+9\left[\left(\sum_{i=2}^nw_i\right)/(n-1)\right]^{1/4} \end{aligned}$  & $\begin{aligned} w_1 & \in [0,1]\\[-.5ex] w_i &=0,  \\[-.5ex] i &= 2,\dots,n \end{aligned}$ & $\begin{aligned} &\text{non-convex,} \\[-.5ex] &\text{non-uniform} \end{aligned}$\\ \bottomrule 

  \end{tabularx}
  }
  \end{table*}

MATLAB's NSGA-II genetic algorithm, a fast sorting and elitist multi-objective genetic algorithm, was used for practical implementation. Initialization was performed through single objective optimization for each objective and random sampling. A population size of 372 was used, with a total of 3,000 maximum generations. The `phenotype' crowding distance metric was used. This setup was validated on the test functions described above. All the problems were benchmarked with $25$ variables ($d = 25$) and two objective functions ($K = 2$), as with the present airfoil optimization problem. The results of our setup on these four benchmark problems are shown in Figure~\ref{fig:benchmark}. The algorithm could accurately capture $ZDT1$, $ZDT2$, and $ZDT4$ and predict $ZDT6$, which is the most complicated due to its non-convex and non-uniform properties, reasonably well.
\begin{figure}[tb!]
  \centering
  \vspace*{1.7em}%
  \includegraphics[width=\linewidth]{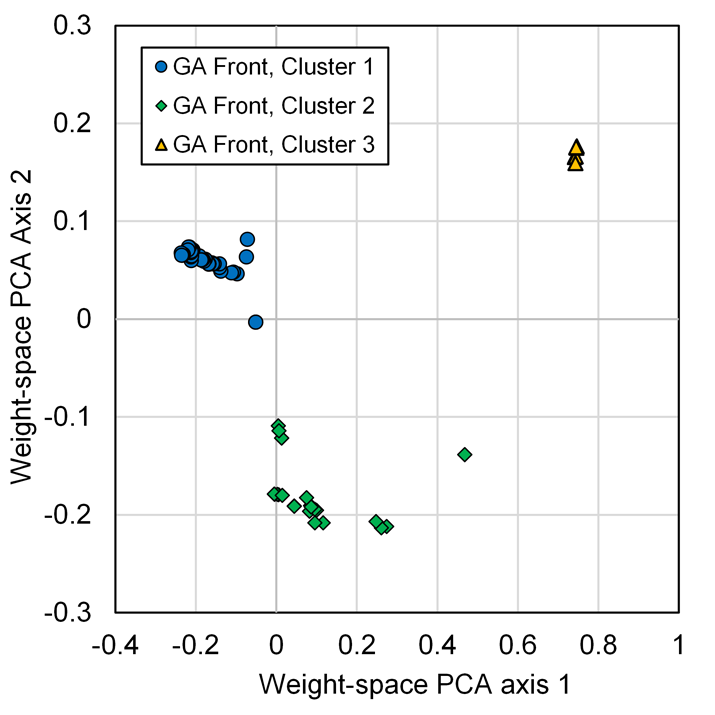}
  \caption{Projection of the 25-dimensional optimal weight vectors to the 2-dimensional subspace spanned by 2 PCA axes of the dominant variance. $k$-means clustering with the cluster size of 3 is used to identify the clusters.}
  \label{fig:pca_cluster}
\end{figure}

\section{Airfoil Shape Clustering} \label{app:pca}
To analyze characteristics of the optimized airfoil shapes in detail, the airfoil shapes on the Pareto-front were classified into three clusters using $k$-means clustering based on the Euclidean distance with $k = 3$. The clustering was performed in the design variable space, or weight space, of $\mathcal{D}^{25}$ rather than in the objective plane because the purpose of clustering was to identify common geometric features over different airfoil shapes as a result of the optimization. The selection of the cluster size was based on the PCA of the optimal weight vector set. 

It should be noted here that the baseline shapes chosen might be linearly dependent. The distances in the PCA weight space, thus, might not be rigorous as a morphed shape on the Pareto-front may be represented by another set of weights. However, this PCA analysis was used only to identify if qualitative classes within the Pareto-front could be found and clustered together and to glean some additional insights of our Pareto-front results.

Figure~\ref{fig:pca_cluster} shows the projection of the 25-dimensional weight vector set to the 2-dimensional subspace spanned by the 2 PCA axes having the first- and second-most variance. The explained variance ratios of PCA axes 1 and 2 are 80.7\% and 14.0\%, respectively. On the other hand, the PCA axis of the third-most variance only accounts for 1.7\% of the variance, affirming that the 2-dimensional projection in Figure~\ref{fig:pca_cluster} adequately scatters the clusters. Based on this observation, $k = 3$ was chosen to be the most appropriate cluster size.

\printbibliography

\end{document}